\documentclass[a4paper,12pt]{article}

\usepackage{listings}
\usepackage{longtable}

\usepackage{graphicx}

\usepackage[font=scriptsize]{caption}

\usepackage{float}

\usepackage{amsmath,amssymb,mathrsfs,eucal}

\usepackage[T1,T2A]{fontenc}

\font\cyrfont=wncyss10

\def\sza{\hbox{\cyrfont X}} 

\newtheorem{thm}{Theorem}

\newtheorem{lem}{Lemma}

\newtheorem{conj}[thm]{Conjecture}

\begin{document}

\title{Behaviour of the order of Tate-Shafarevich groups for the quadratic 
twists of elliptic curves}

\author{Andrzej D\k{a}browski and Lucjan Szymaszkiewicz} 

\date{}

\maketitle{}

{\it Abstract}. We present the results of our search for the 
orders of Tate-Shafarevich groups for the quadratic twists of elliptic 
curves. We formulate a general conjecture, giving for a fixed elliptic 
curve $E$ over $\Bbb Q$ and positive integer $k$, an asymptotic formula 
for the number of quadratic twists $E_d$, $d$ positive square-free integers 
less than $X$, with finite group $E_d(\Bbb Q)$ and $|\sza(E_d(\Bbb Q))| = k^2$. 
This paper continues the authors previous investigations concerning orders 
of Tate-Shafarevich groups in quadratic twists of the curve $X_0(49)$. 
In section 8 we exhibit $88$ examples of rank zero elliptic curves 
with $|\sza(E)| > 63408^2$, which was the largest previously known 
value for any explicit curve. Our record is an elliptic curve $E$ with 
$|\sza(E)| = 1029212^2$. 

\bigskip 
Key words: elliptic curves, Tate-Shafarevich group, 
Cohen-Lenstra heuristics, distribution 
of central $L$-values

\bigskip 
2010 Mathematics Subject Classification: 11G05, 11G40, 11Y50

\section{Introduction}

Let $E$ be an elliptic curve defined over $\Bbb Q$ of conductor $N_E$, 
and let $L(E,s)$ denote its $L$-series. 
Let $\sza(E)$ be the Tate-Shafarevich group of $E$, 
$E(\Bbb Q)$ the group of rational points, and $R(E)$ the regulator,  
with respect to the N\'eron-Tate height pairing. 
Finally, let $\Omega_E$ be the least positive real period of the N\'eron 
differential of a global minimal Weierstrass equation for 
$E$, and define $C_{\infty}(E)=\Omega_E$ or $2\Omega_E$ 
according as $E(\Bbb R)$ is connected or not, and let $C_{\text{fin}}(E)$ 
denote the product of the Tamagawa factors of $E$ at the bad primes. 
The Euler product defining $L(E,s)$ converges for $\text{Re}\,s>3/2$. 
The modularity conjecture, proven by Wiles-Taylor-Diamond-Breuil-Conrad,   
implies that $L(E,s)$ has an analytic continuation to an entire function. 
The Birch and Swinnerton-Dyer conjecture relates the arithmetic data 
of $E$ to the behaviour of $L(E,s)$ at $s=1$. 

\begin{conj} (Birch and Swinnerton-Dyer) (i) $L$-function $L(E,s)$ has a 
zero of order $r=\text{rank}\;E(\Bbb Q)$ at $s=1$, 

(ii) $\sza(E)$ is finite, and 
$$ 
\lim_{s\to1}\frac{L(E,s)}{(s-1)^r} =
\frac{C_{\infty}\, (E)C_{\text{fin}}(E)\, R(E)\, |\sza(E)|}{|E(\Bbb Q)_{\text{tors}}|^2}. 
$$
\end{conj} 
If $\sza(E)$ is finite, the work of Cassels and Tate shows that its order 
must be a square. 

The first general result in the direction of this conjecture was proven for 
elliptic curves $E$ with complex multiplication by Coates and Wiles in 1976 
\cite{CW}, who showed that if $L(E,1)\not =0$, then the group $E(\Bbb Q)$ is finite. 
Gross and Zagier \cite{GZ} showed that if $L(E,s)$ has a first-order zero at $s=1$, 
then $E$ has a rational point of infinite order. Rubin \cite{Rub} proves that 
if $E$ has complex multiplication and $L(E,1)\not =0$, then $\sza(E)$ is finite. 
Let $g_E$ be the rank of $E(\Bbb Q)$ and let $r_E$ the order of the zero 
of $L(E,s)$ at $s=1$. Then Kolyvagin \cite{Kol} proved that, if $r_E\leq 1$, 
then $r_E=g_E$ and $\sza(E)$ is finite.  Very recently, Bhargava, Skinner 
and Zhang \cite{BhSkZ} proved that at least $66.48 \%$ of all elliptic curves 
over $\Bbb Q$, when ordered by height, satisfy the weak form of the Birch 
and Swinnerton-Dyer conjecture, and have finite Tate-Shafarevich group.

When $E$ has complex multiplication by the ring of integers of an imaginary 
quadratic field $K$ and $L(E,1)$ is non-zero, the $p$-part of the Birch and 
Swinnerton-Dyer conjecture has been established by Rubin \cite{Rub2} for 
all primes $p$ which do not divide the order of the group of roots of unity 
of $K$. Coates et al. \cite{cltz} \cite{coa}, and Gonzalez-Avil\'es 
\cite{G-A} showed that there is 
a large class of explicit quadratic twists of $X_0(49)$ whose 
complex $L$-series does not vanish at $s=1$, and for which the full 
Birch and Swinnerton-Dyer conjecture is valid (covering the case $p=2$ 
when $K=\Bbb Q(\sqrt{-7})$). 
The deep results by Skinner-Urban (\cite{SkUr}, Theorem 2) (see also Theorem 7 
in section 8.4 below) 
allow, in specific cases (still assuming $L(E,1)$ is non-zero), to establish $p$-part of 
the Birch and Swinnerton-Dyer conjecture for elliptic curves  without complex multiplication 
for all odd primes $p$ (see examples in section 8.4 below, and section 3 in \cite{DSz}).

The numerical studies and conjectures by 
Conrey-Keating-Rubinstein-Snaith \cite{CKRS}, 
Delaunay \cite{Del0}\cite{Del}, Watkins \cite{Wat}, 
Radziwi\l\l-Soun\-dararajan \cite{RS}  
(see also the papers \cite{dw}\cite{DSz}\cite{djs}, and references therein) 
substantially extend the systematic tables given by Cremona. 

This paper continues the authors previous investigations concerning 
orders of Tate-Shafarevich groups in quadratic twists of the 
curve $X_0(49)$. We present the results of our search for the orders 
of Tate-Shafarevich groups for additional four elliptic curves (two 
with CM, and two without CM), for the same large ranges of the index 
(namely, $32\cdot 10^9$). Our results support a general conjecture 
(Conjecture 2), giving for a fixed elliptic 
curve $E$ over $\Bbb Q$ and positive integer $k$, an asymptotic formula 
for the number of quadratic twists $E_d$, $d$ positive square-free integers 
less than $X$, with finite group $E_d(\Bbb Q)$ and $|\sza(E_d(\Bbb Q))| = k^2$. 
Let us formulate explicitly the conjecture.  Let $f_E(k,X)$ denote the number 
of positive square-free integers $d\leq X$, such that $(d,N_E)=1$, 
$L(E_d,1)\not=0$, and $|\sza(E_d)|=k^2$. 

\bigskip 
\noindent 
{\bf Conjecture} {\it  Let $E$ be an elliptic curve defined over $\Bbb Q$. 
For any positive integer $k$ there are constants $c_k(E)\geq 0$ and 
$d_k(E)$ such that}  
$$
f_E(k,X) \sim c_k(E) X^{3/4}(\log X)^{d_k(E)}, \quad \text{as} 
\quad X\to\infty.  
$$ 

The main results supporting the conjecture are reported in section 3;  
additional support  is also given in subsection 7.1. 
In section 9 we give examples of elliptic curves $E$ with $c_r(E)=0$, 
for some values $r$. 

All the experiments concerning statistics of the $L$-values of quadratic 
twists of $X_0(49)$ (and related orders of Tate-Shafarevich groups) done in 
\cite{djs}, are also confirmed for these four elliptic curves 
(see sections 4 - 6).  

It has long been known that the order of $\sza(E)[p]$ can be arbitrarily 
large for elliptic curves $E$ defined over $\Bbb Q$ and $p=2, 3$ 
(for $p=3$, the result 
is due to J. Cassels (\cite{Cas}), and for $p=2$ it is due to F. McGuinness 
(\cite{McG}), but no similar result is known for $p > 3$. 
We also stress that it has not yet been proven that there exist elliptic 
curves $E$ defined over $\Bbb Q$ for which $\sza(E)[p]$ is non-zero for 
arbitrarily large primes $p$. In our earlier papers, we have investigated 
(see \cite{dw}, \cite{djs}) 
some numerical examples of $E$ defined over $\Bbb Q$ for which $L(E, 1)$ 
is non-zero and the order of $\sza(E)$ is large. We extend these numerical 
results here in section 8, with the largest proved examples of $\sza(E)$ 
having order $1029212^2 = 2^4\cdot 79^2 \cdot 3257^2$.  


\bigskip 
We heartily thank John Coates for his remarks, suggestions and 
many corrections.

\bigskip 

This research was supported in part by PL-Grid Infrastructure.
Our computations were carried out in 2016 on the Prometheus
supercomputer via PL-Grid infrastructure.  For the calculations 
in section 8 we also used the HPC cluster HAL9000 and descop 
computers Core(TM) 2 Quad Q8300 4GB/8GB, all located at the 
Department of Mathematics and Physics of Szczecin University.

\section{Formulae for the orders of $\sza(E_d)$ when $L(E_d,1)\not=0$}

We can compute $L(E_d,1)$ when it is non-zero for a huge range of 
positive square-free integers $d$ thanks to the 
remarkable ideas discovered by Waldspurger, and worked out explicitly 
in particular cases by many authors. These ideas show that 
$L(E_d,1)$, when it is non-zero, is essentially equal to the 
$d$-th Fourier coefficient of an explicit modular form of weight $3/2$.  
We now recall some details for four elliptic curves (named $A$, $B$, 
$C$, and $D$ below).  

\bigskip 

{\it Notation}. Let $q:=e^{2\pi iz}$, 
$\Theta(z)=\sum_{n\in\Bbb Z}q^{n^2}$, 
$\eta(z)=e^{\pi iz/12}\prod_{n=1}^{\infty}(1-e^{2\pi i nz})$.

\bigskip 

{\it Example 1}. \cite{Tu} Let $A: y^2=x^3-x$; for any square-free integer 
$d$ consider the quadratic twist $A_d: y^2=x^3-d^2x$. 

Let $f_i(z):=\eta(8z)\eta(16z)\Theta(2^iz) = \sum a_i(n)q^n$ ($i=1,2$). 
Consequently, for odd $d$ we have  
$
a_i(d)=|\{(x,y,z)\in\Bbb Z^3: d=2ix^2+y^2+32z^2\}| - {1\over 2}
|\{(x,y,z)\in\Bbb Z^3: d=2ix^2+y^2+8z^2\}|. 
$
If $d\geq 1$ is an odd square-free integer, then Tunnell proved that
$$
L(A_{id},1)={2^{i-1}a_i(d)^2\Omega_A\over 4\sqrt{2^{i-1}d}}. 
$$
Therefore assuming the Birch and Swinnerton-Dyer conjecture, $A_{id}$ has 
$\Bbb Q$-rank zero iff $a_i(d)\not=0$. In addition, if $a_i(d)\not=0$, then
$$
|\sza(A_{id})|=\left({a_i(d)\over \tau(d)}\right)^2,
$$
where $\tau(d)$ denotes the number of divisors of $d$. Let $a(d):=a_1(d)$.

{\it Example 2.} \cite{Fr2},\cite{Fe} Let $B: y^2=x^3-1$. For any square-free 
positive integer $d\equiv 1(\text{mod}\;4)$, $(d,6)=1$, consider the quadratic 
twist $B_d: y^2=x^3-d^3$. 
Let $a(d)$ denote the $d$-th Fourier coefficient of $\eta^2(12z)\Theta(z)$. 
Frey (\cite{Fr2}, page 232) proves the following result 
$$
L(B_d,1)=
\begin{cases} 
0
&\text{if $d\equiv 3(\text{mod}\;4)$}\\
a(d)^2{L(B,1)\over\sqrt{d}}
&\text{if $d\equiv 1(\text{mod}\;24)$}\\
\left({a(d)\over a(13)}\right)^2\sqrt{{13\over d}}L(B_{13},1)
&\text{if $d\equiv 13(\text{mod}\;24)$}\\
\left({a(d)\over a(5)}\right)^2\sqrt{{5\over d}}L(B_5,1)
&\text{if $d\equiv 5(\text{mod}\;24)$}\\
\left({a(d)\over a(17)}\right)^2\sqrt{{17\over d}}L(B_{17},1)
&\text{if $d\equiv 17(\text{mod}\;24)$}
\end{cases}  
$$ 

One can show \cite{Fe} that $a(d)={1\over 2}\sum(-1)^n$, where the sum is taken over all 
$m,n,k\in\Bbb Z$ satisfying $m^2+n^2+k^2=d$, $3\not|m$, $3|n$ and $2\not|m+n$.

Let 
$l=l_1+[{l_2+1\over 2}]$, where  $l_i=|\{p|d: p\equiv i(\text{mod}\;3)\}|$. 
Assuming the Birch and Swinnerton-Dyer conjecture, we obtain
$$
|\sza(B_d)|=\left({a(d)\over 2^l}\right)^2,
$$
if $a(d)\not=0$.

{\it Example 3.} \cite{Chen} Let $C: y^2=x^3+4x^2-16$ be an elliptic curve of conductor $176$. 
Consider the family $C_d: y^2=x^3+4dx^2-16d^3$, 
where $d$ runs over positive, odd, square-free integers, satisfying 
$({d\over 11})=1$, i.e. such that $d\equiv 1,3,4,5,9 (\text{mod}\, 11)$.  
Let $a(d)={n_d-m_d\over 2}$, where 
$
n_d=|\{(x,y,z)\in \Bbb Z^3: d=x^2+11y^2+11z^2\}|, 
$ 
and 
$
m_d=|\{(x,y,z)\in \Bbb Z^3: d=3x^2+2xy+4y^2+11z^2\}|.  
$ 
All the groups $C_d(\Bbb Q)_{tors}$ are trivial, hence 
assuming the Birch and Swinnerton-Dyer conjecture, we obtain 
$$
|\sza(C_d)|={a(d)^2 \over C_{fin}(C_d)},
$$
if $a(d)\not=0$.  Moreover, in this case we have 
$$
L(C_d,1)={2^{s(d)}a(d)^2\Omega_C\over \sqrt{d}}, 
$$
where $s(d)=0$ if $d\equiv 1 (\text{mod}\, 4)$ and 
$s(d)=1$ if $d\equiv 3 (\text{mod}\, 4)$.

{\it Example 4.} \cite{Chen} Let $D: y^2=x^3-x^2-8x-16$ be an elliptic curve of conductor $112$. 
Consider the family $D_d: y^2=x^3-dx^2-8d^2x-16d^3$, 
where $d$ runs over positive, odd, square-free integers, satisfying 
$({d\over 7})=1$, i.e. such that $d\equiv 1,2,4 (\text{mod}\, 7)$. 
Let $a(d)={n_d-m_d\over 2}$, where 
$
n_d=|\{(x,y,z)\in \Bbb Z^3: d=x^2+14y^2+14z^2\}|, 
$ 
and 
$
m_d=|\{(x,y,z)\in \Bbb Z^3: d=2x^2+7y^2+14z^2\}|.  
$ 
Let us stress, that in case $d\equiv 3 (\text{mod}\, 8)$ we necessarily 
have $a(d)=0$. Let $s(d)=2$ if $d\equiv 1 (\text{mod}\, 4)$ and 
$s(d)=1$ if $d\equiv 3 (\text{mod}\, 4)$. 
All the groups $D_d(\Bbb Q)_{tors}$ have order two, hence 
assuming the Birch and Swinnerton-Dyer conjecture, we obtain 
$$
|\sza(D_d)|={2^{s(d)}a(d)^2 \over C_{fin}(D_d)}, 
$$
if $a(d)\not=0$. Moreover, in this case we have 
$$
L(D_d,1)={a(d)^2\Omega_D\over \sqrt{d}}.  
$$

\noindent 
{\it Definitions}. We say, that a positive square-free odd integer 
$d$ satisfies:  

\bigskip 

(i) condition ($*_A$), if $d\equiv 1,3 (\text{mod}\, 8)$; 

(ii) condition ($*_B$), if $d\equiv 1,5 (\text{mod}\, 12)$; 

(iii) condition ($*_C$), if $d\equiv 1,3,4,5,9 (\text{mod}\, 11)$; 

(iv) condition ($*_D$), if $d\equiv 1,2,4 (\text{mod}\, 7)$ and 
$d\not\equiv 3 (\text{mod}\, 8)$;   

(v) condition ($**_E$), if it satisfies condition ($*_E$), and $a(d)\not=0$.

\section{Frequency of orders of $\sza$}

Our data contains values of $|\sza(E_d)|$ for $E\in\{A,B,C,D\}$, 
and $d\leq 32\cdot 10^9$ satisfying ($**_E$). 
The data involves the proven odd orders of $|\sza(A_d)|$, and the 
proven (prime to $6$) orders of $|\sza(B_d)|$. 
The non-trivial values 
of $|\sza(C_d)|$ and $|\sza(D_d)|$ are (mostly) the conjectural ones. 
Let $k_E$ denote 
the largest positive integer such that for all $k\leq k_E$ 
there is quadratic twist $E_{d_k}$ (with $d_k$ as above) with 
$|\sza(E_{d_k})|=k^2$. Let $K_E$ denote the largest positive 
integer $k$ such that for some $d_k$ as above, we have 
$|\sza(E_{d_k})|=k^2$. 

\begin{center}
\footnotesize
\begin{longtable}{|r|r|r|r|r|r|} 
\hline 
\multicolumn{1}{|c|}{$E$} 
& 
\multicolumn{1}{|c|}{$k_E$} 
& 
\multicolumn{1}{|c|}{$K_E$}  
\\ 
\hline 
\endhead 
\hline 
\hline
\endfoot

\hline 
\hline
\endlastfoot

$A$ & 2277 & 2783   \\
$B$ & 2037 & 3571   \\ 
$C$ & 1727 & 4235   \\ 
$D$ & 1914 & 2667   \\ 

\end{longtable}
\end{center} 

Note that $3571$ is a prime. From our data it follows that 
$|\sza(C_{26650821201})|=3917^2$, with $3917$ the largest known 
(at the moment) prime dividing the order of $\sza(E_d)$ of 
an elliptic $E\in\{A, B, C, D, X_0(49)\}$.

Our calculations strongly suggest that 
for any positive integer $k$ there are infinitely many positive 
integers $d$ satisfying condition 
($*_E$), such that $E_d$ has rank zero and $|\sza(E_d)| = k^2$. 
Below we will state a more precise conjecture.

Let $f_E(X)$ denote the number of integers 
$d\leq X$, satisfying ($**_E$) and 
such that $|\sza(E_d)| = 1$. Let 
$g_E(X)$ denote the number of integers $d\leq X$, 
satisfying ($*_E$) and such that $L(E_d,1) = 0$. 

We obtain the following graphs of the functions $f_E(X)/g_E(X)$, for 
$E\in\{A, B, C, D\}$.

\begin{figure}[H]  
\centering 
\includegraphics[trim = 0mm 20mm 0mm 15mm, clip, scale=0.4]{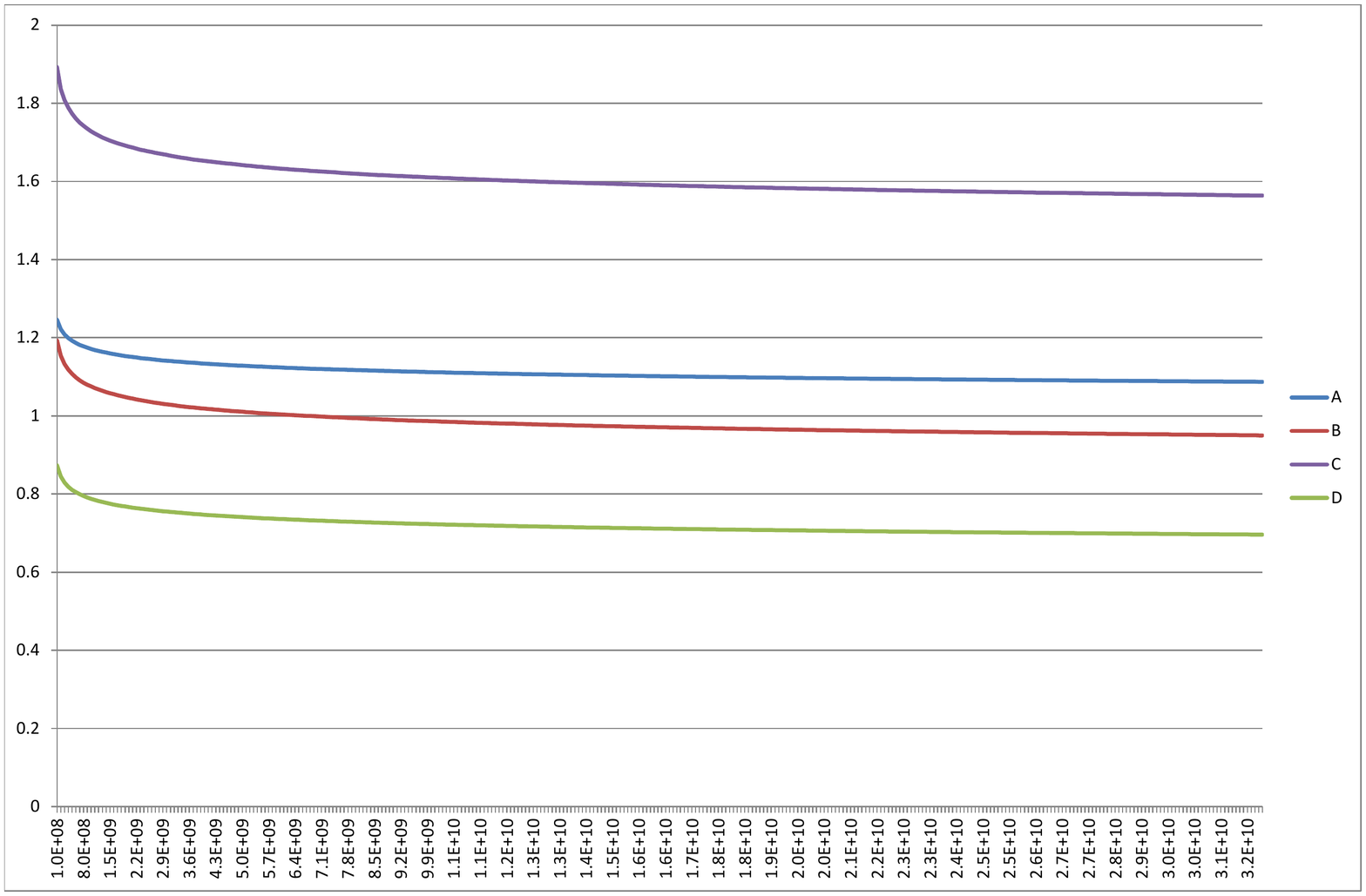}
\caption{Graphs of the functions $f_E(X)/g_E(X)$, for $E\in\{A,B,C,D\}$.} 
\end{figure}

We expect that $f_E(X)/g_E(X)$ tend to constants dependent on $E$ 
(the constant is $1$ for $E=X_0(49)$, see section 11 in \cite{djs}).

We expect (Delaunay-Watkins \cite{DW}, Heuristics 1.1): 
$$
g_E(X) \sim c_E X^{3/4}(\log X)^{b_E+{3\over 8}}, \quad \text{as} 
\quad X\to\infty, 
$$
where $c_E>0$, and there are four different possibilities for $b_E$, 
largery dependent on the rational $2$-torsion structure of $E$. 
Hence, we may expect similar asymptotic formula for $f_E(X)$ as well.

Now let $f_E(k,X)$ denote the number of integers 
$d\leq X$, satisfying ($**_E$) and 
such that $|\sza(E_d)| = k^2$. 
Let $F_E(k,X):={f_E(X)\over f_E(k,X)}$. We obtain the following graphs 
of the functions $F_E(k,X)$ for $E\in\{A, B, C, D\}$ and $k=2,3,4,5,6,7$.

\begin{figure}[H]  
\centering
\includegraphics[trim = 0mm 20mm 0mm 15mm, clip, scale=0.4]{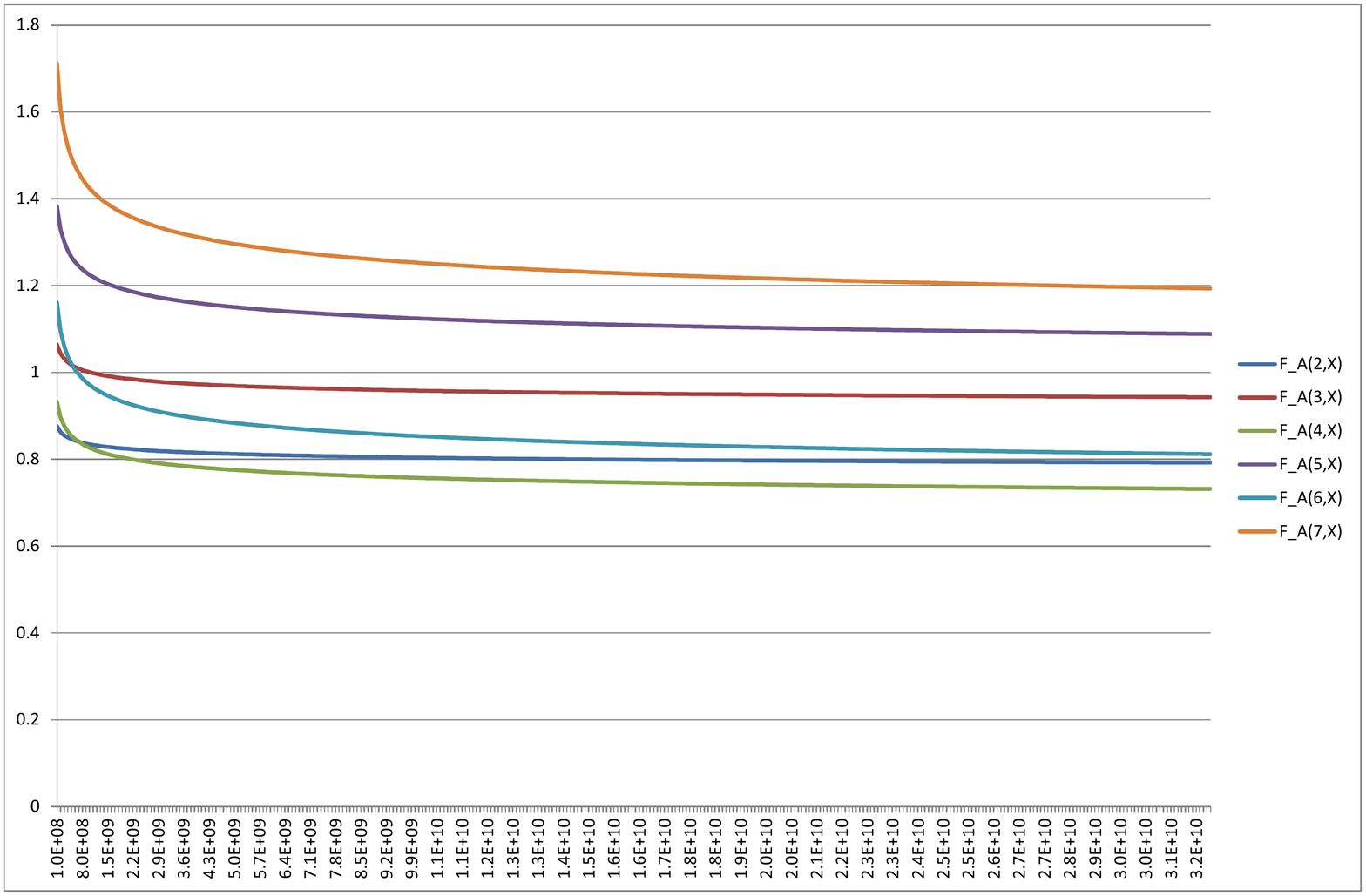}
\caption{Graphs of the functions $F_A(k,X)$ for $k=2,3,4,5,6,7$.} 
\end{figure}

\begin{figure}[H]  
\centering
\includegraphics[trim = 0mm 20mm 0mm 15mm, clip, scale=0.4]{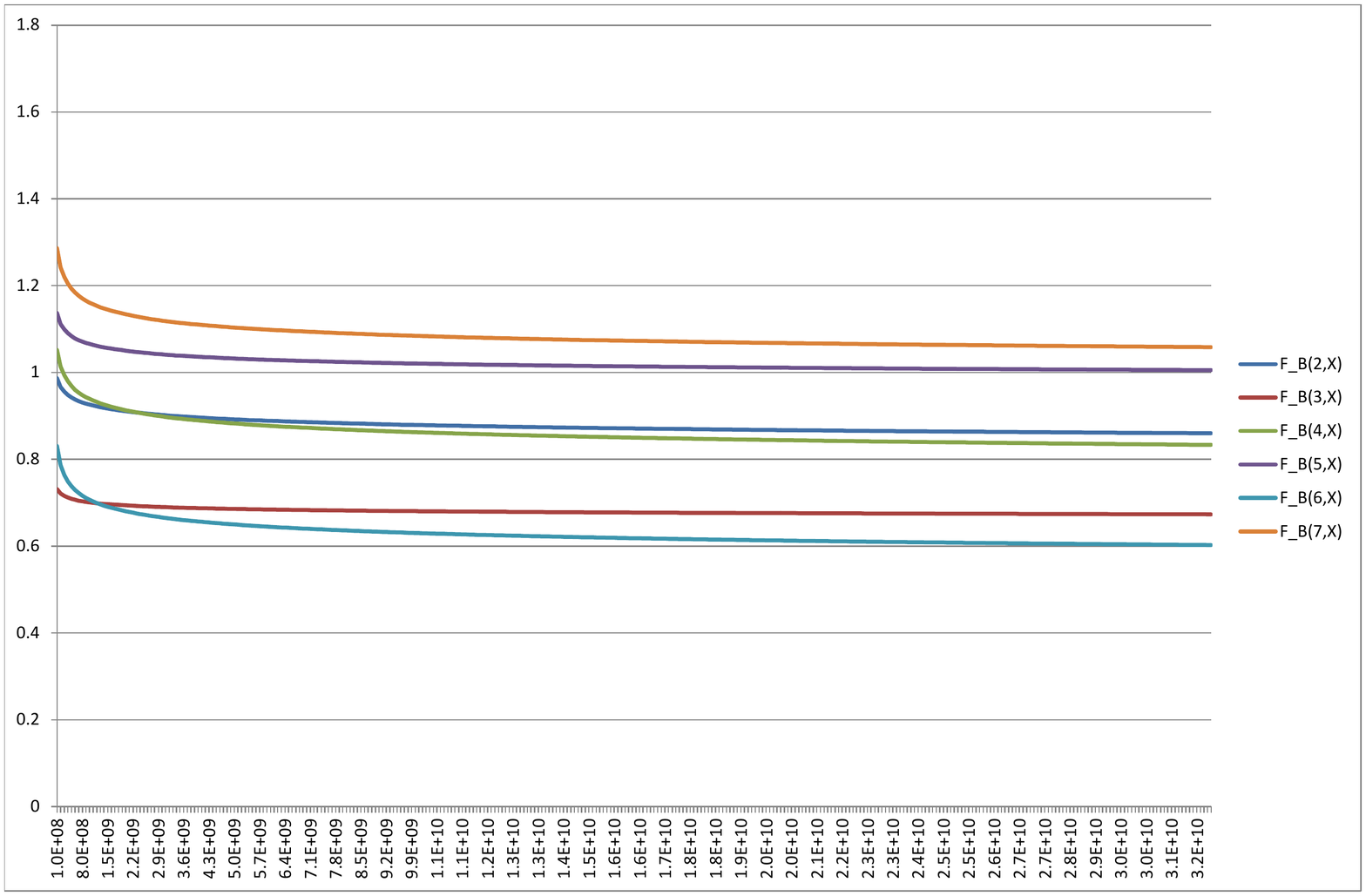}
\caption{Graphs of the functions $F_B(k,X)$ for $k=2,3,4,5,6,7$.} 
\end{figure}

\begin{figure}[H]  
\centering
\includegraphics[trim = 0mm 20mm 0mm 15mm, clip, scale=0.4]{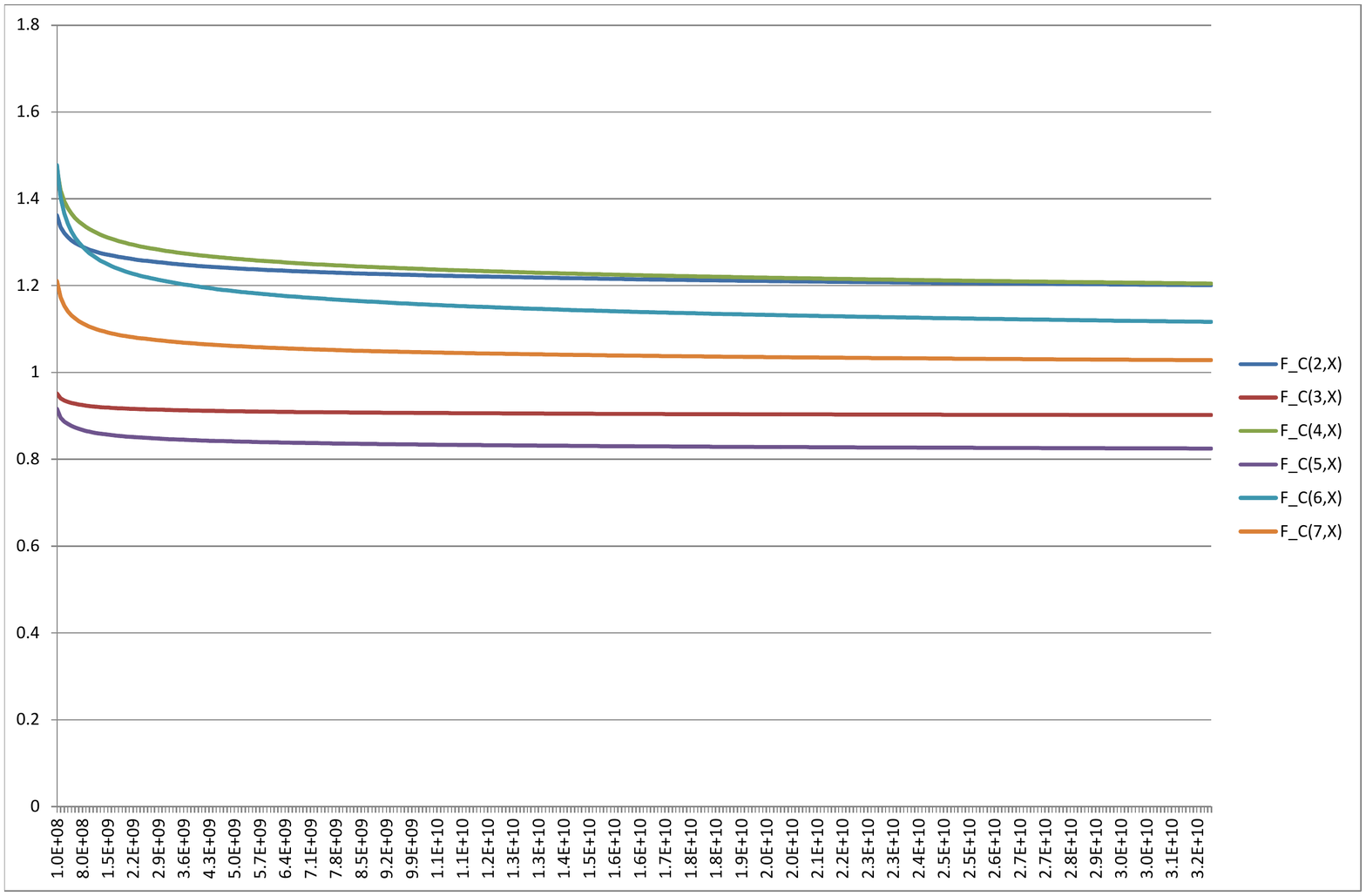}
\caption{Graphs of the functions $F_C(k,X)$ for $k=2,3,4,5,6,7$.} 
\end{figure} 

\begin{figure}[H]  
\centering
\includegraphics[trim = 0mm 20mm 0mm 15mm, clip, scale=0.4]{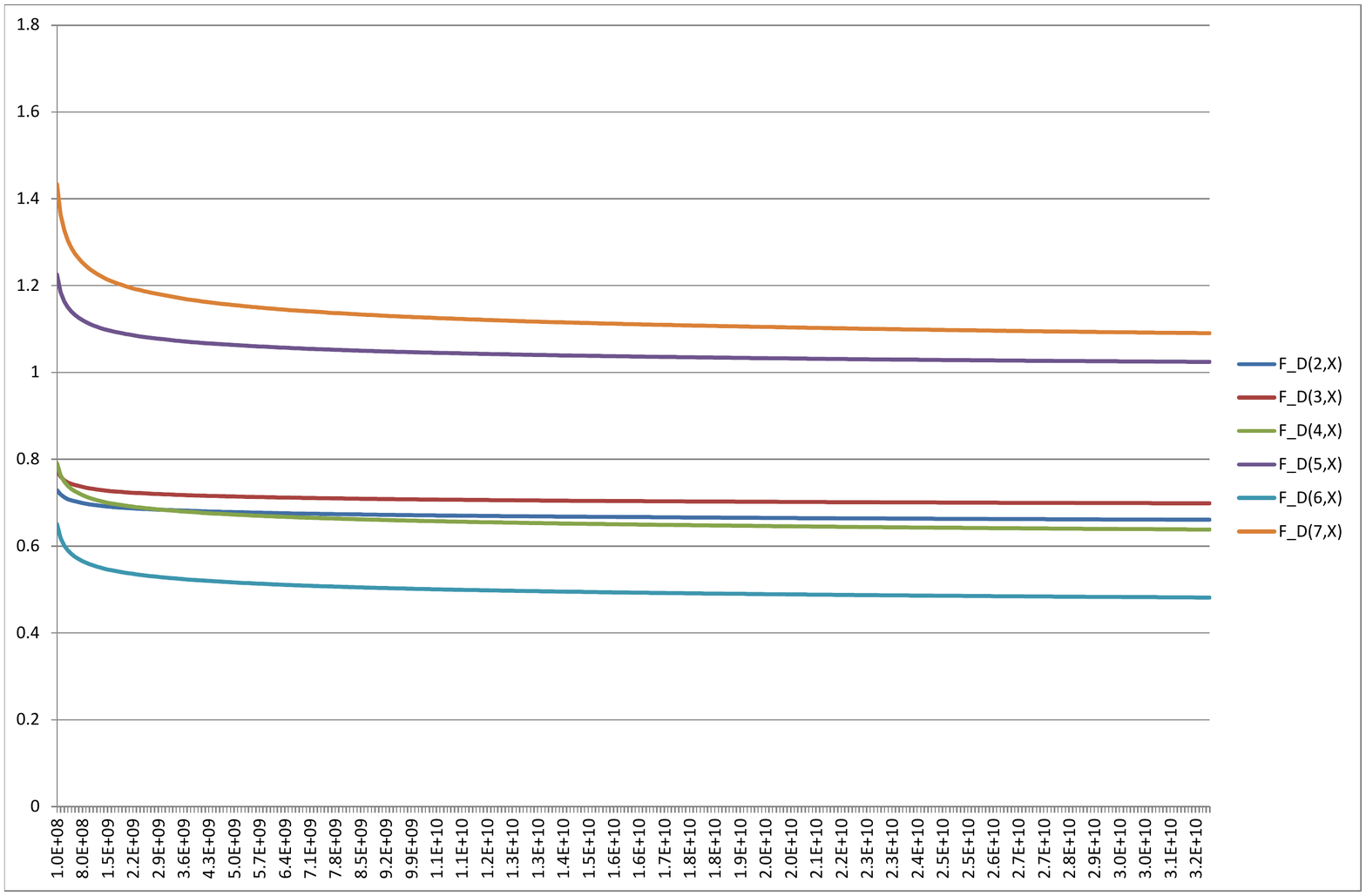}
\caption{Graphs of the functions $F_D(k,X)$ for $k=2,3,4,5,6,7$.} 
\end{figure}

The above calculations suggest the following general conjecture 
(compare Conjecture 8 in \cite{djs} for the case of quadratic 
twists of the curve $X_0(49)$). Let $f_E(k,X)$ denote the number 
of positive square-free integers $d\leq X$, such that $(d,N_E)=1$, 
$L(E_d,1)\not=0$, and $|\sza(E_d)|=k^2$.

\begin{conj} Let $E$ be an elliptic curve defined over $\Bbb Q$. 
For any positive integer $k$ there are constants $c_k(E)\geq 0$ and 
$d_k(E)$ such that 
$$
f_E(k,X) \sim c_k(E) X^{3/4}(\log X)^{d_k(E)}, \quad \text{as} 
\quad X\to\infty. 
$$
\end{conj} 
In section 9 we give examples of elliptic curves $E$ with 
$c_{2m+1}(E)=0$ and with $c_{4m+2}(E)=0$. 

Note that Park, Poonen, Voight, and Wood (\cite{PPVW}, (11.2.2)) have 
formulated an analogous conjecture for the family of all elliptic curves 
over the rationals, ordered by height.

\section{Cohen-Lenstra heuristics for the order of $\sza$} 

Delaunay \cite{Del} has considered Cohen-Lenstra heuristics for the order of 
Tate-Shafarevich group. He predicts, among others, that in the rank zero case, 
the probability that $|\sza(E)|$ of a given elliptic curve $E$ over $\Bbb Q$ is 
divisible by a prime $p$ should be 
$
f_0(p):=1-\prod_{j=1}^{\infty}(1-p^{1-2j})={1\over p}+{1\over p^3}+... \,.    
$
Hence, $f_0(2)\approx 0.580577$, $f_0(3)\approx 0.360995$, $f_0(5)\approx 0.206660$, 
$f_0(7)\approx 0.145408$,  $f_0(11)\approx 0.092$, 
and so on.

Let $F_E(X)$ denote the number of positive integers 
$d\leq X$,  
satisfying ($*_E$) and such that $L(E_d,1)\not=0$. 
Let $F_E(p,X)$ denote the number of positive integers 
$d\leq X$, satisfying ($**_E$), such that  
$|\sza(E_d)|$ is divisible by $p$. 
Let $f_E(p,X):={F_E(p,X)\over F_E(X)}$, and  
$f_E(p):=f_E(p,32\cdot 10^9)$. We obtain the following table

\begin{center}
\footnotesize
\begin{longtable}{|r|r|r|r|r|r|r|r|} 
\hline 
\multicolumn{1}{|c|}{$E$} 
& 
\multicolumn{1}{|c|}{$f_E(2)$} 
& 
\multicolumn{1}{|c|}{$f_E(3)$} 
& 
\multicolumn{1}{|c|}{$f_E(5)$} 
& 
\multicolumn{1}{|c|}{$f_E(7)$} 
& 
\multicolumn{1}{|c|}{$f_E(11)$} 
\\ 
\hline 
\endhead 
\hline 
\multicolumn{5}{|r|}{{Continued on next page}} \\
\hline
\endfoot

\hline 
\hline
\endlastfoot

$A$ & 0.565173 & 0.348417 &  0.192196 &  0.130318 &  0.076544 \\
$B$ & 0.500459 & 0.427769 &  0.197255 &  0.135517 &  0.081501 \\ 
$C$ & 0.387009 & 0.355532 &  0.233508 &  0.139171 &  0.085150 \\ 
$D$ & 0.607500 & 0.424331 &  0.193023 &  0.131217 &  0.077425 \\ 

\end{longtable}
\end{center}

The papers of Quattrini \cite{Qu1}\cite{Qu2} make a correction to 
Delaunay's heuristics for $p$-divisibility of $|\sza(E_d)|$ in the 
family of quadratic twists of a given elliptic curve $E$ of square-free 
conductor for odd primes dividing the order of $E(\Bbb Q)_{tors}$. 
The author gives an explanation of why and when the original Cohen-Lenstra 
heuristics should be used for the prediction of the $p$-divisibility of 
the order of $\sza(E_d)$. Roughly speaking, the proportion of 
values of $|\sza(E_d)|$ divisible by a prime number $p$ among 
(imaginary) quadratic twists of $E$ is significantly bigger when $E$ 
has a $\Bbb Q$-rational point of order $p$, than in the general case 
where $|E(\Bbb Q)_{tors}|$ is not divisible by $p$. 
In our situation, $|B_{-1}(\Bbb Q)_{tors}|=|D_{-1}(\Bbb Q)_{tors}|=6$, 
$|C_{-1}(\Bbb Q)_{tors}|=5$, the original Cohen-Lenstra predictions for 
$p=3$ (resp. for $p=5$ are $\approx 0.439$ (resp. $\approx 0.239$), 
and it explains why the values $f_B(3)$, 
$f_D(3)$, and $f_C(5)$ deviate from Delaunay's predictions. We have no 
explanation why the values  $f_B(2)$, $f_C(2)$, $f_D(2)$ deviate from 
the expected one.

\section{Numerical evidence for Delaunay asymptotic formulae}

Let 
$M_E(T):={1\over T^*}\sum |\sza(E_d)|$, 
where the sum is over positive integers $d\leq T$, 
satisfying ($**_E$), and $T^*$ denotes the number of 
terms in the sum. Delaunay (\cite{Del0}, Conjecture 6.1) has conjectured 
that 
$$
M_E(T) \sim c_E T^{1/2}(\log T)^{t_E}, \quad \text{as} 
\quad T\to\infty,  
$$ 
where $-1 < t_E < 1$ largery depend on the rational $2$-structure of $E$ 
(for instance, $t_E=-5/8$ for $E=B$ or $D$). 
If we restrict to prime twists, then we obtain a similar conjecture, 
but without the $\log$ term (Conjecture 4.2 in \cite{Del0}). 

Let $N_E(T)$ be a subsum of $M_E(T)$, restricted to prime twists. 
Let $f_E(T):={M_E(T)\over T^{1/2}}$, and 
$g_E(T):={N_E(T)\over T^{1/2}}$. We obtain the following pictures confirming 
the conjectures 4.2 and 6.1 in \cite{Del0} for the curves $A$, $B$, $C$, 
and $D$ (compare numerical evidence for the curve $X_0(49)$ done in \cite{djs}).

\begin{figure}[H]  
\centering
\includegraphics[trim = 0mm 20mm 0mm 15mm, clip, scale=0.4]{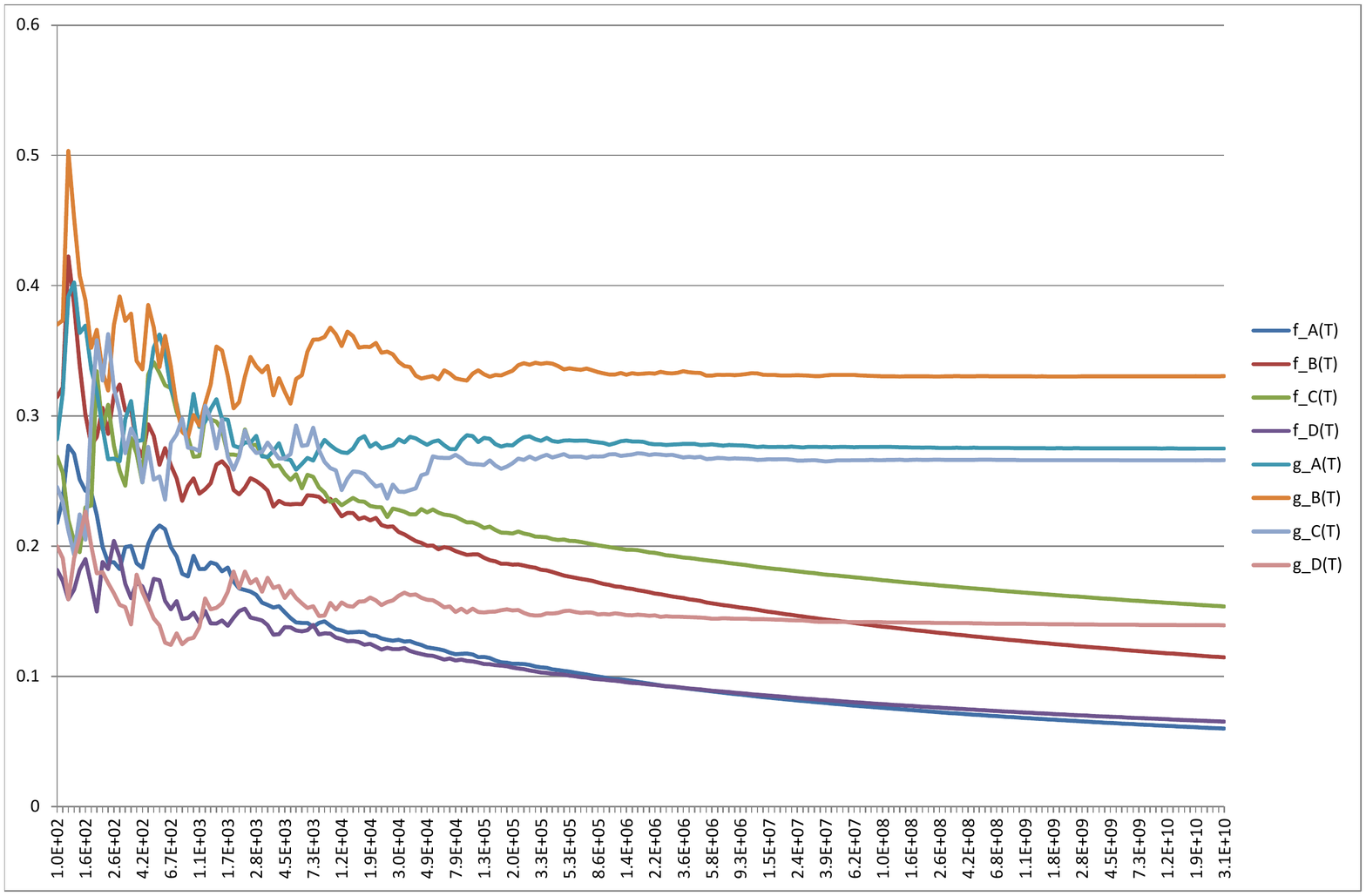} 
\caption{Graphs of the functions $f_E(T)$ and $g_E(T)$, $E\in\{A,B,C,D\}$,  
using the geometric sequence of arguments.} 
\end{figure} 

Let us modify the functions $f_E(T)$ for the curves $B$ and $D$, 
using the logarithmic factor predicted by Delaunay conjecture:  
$f^*_E(T):={(\log T)^{5/8}M_E(T)\over T^{1/2}}$.

\begin{figure}[H]  
\centering
\includegraphics[trim = 0mm 20mm 0mm 15mm, clip, scale=0.4]{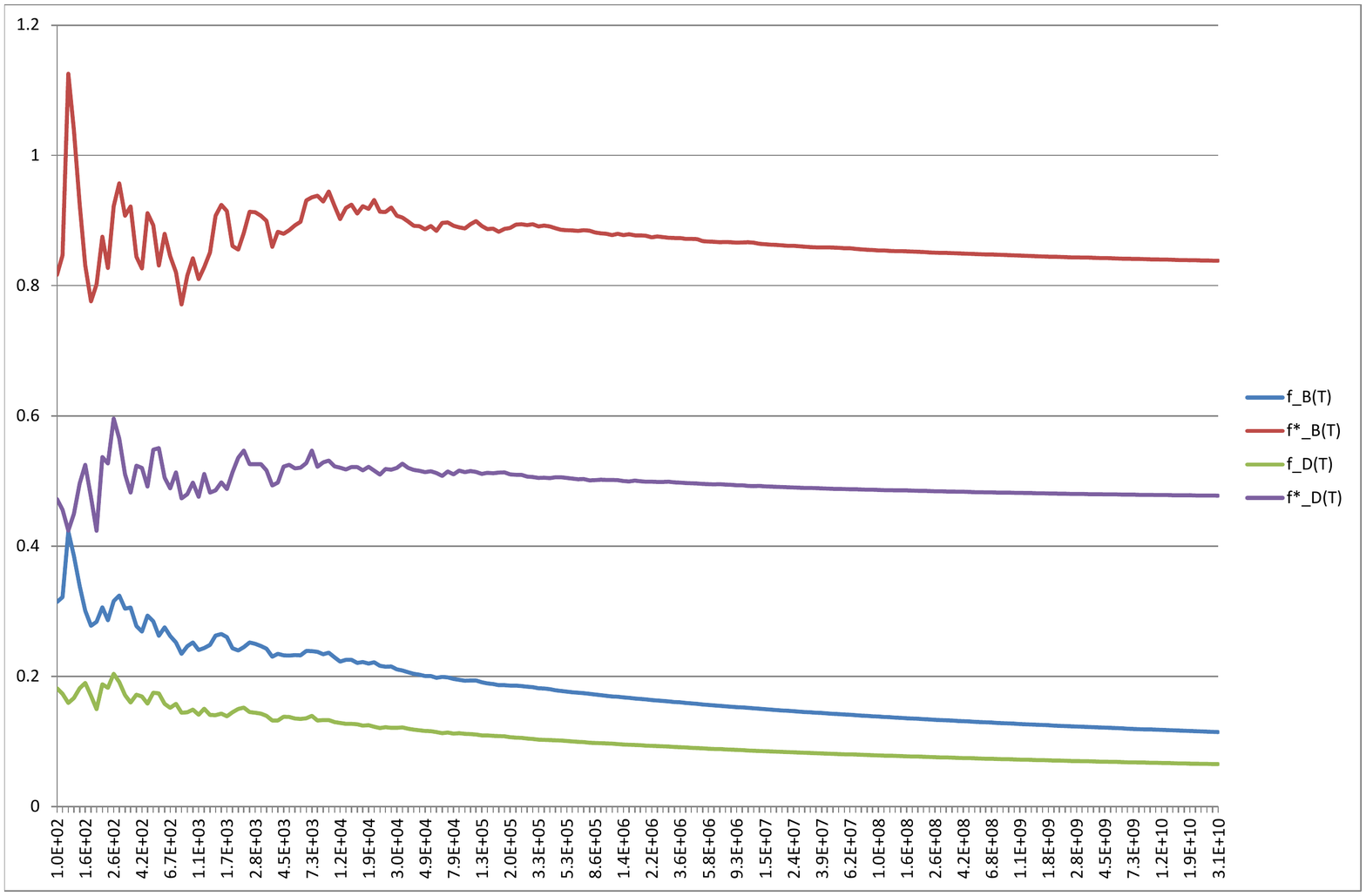} 
\caption{Graphs of the functions $f^*_E(T)$ and $f_E(T)$, $E\in\{B,D\}$,  
using the geometric sequence of arguments.} 
\end{figure}

\section{Distributions of $L(E_d,1)$ and $|\sza(E_d)|$}

\subsection{Distribution of $L(E_d,1)$} 

It is a classical result (due to Selberg) that the values of 
$\log |\zeta({1\over 2}+it)|$ follow a normal distribution.

Let $E$ be any elliptic curve defined over $\Bbb Q$. 
Let $\cal E$ denote the set of all fundamental discriminants $d$ with $(d,2N_E)=1$ and 
$\epsilon_E(d)=\epsilon_E\chi_d(-N_E)=1$, where $\epsilon_E$ is the root number of $E$ 
and $\chi_d=(d/\cdot)$.  Keating and Snaith \cite{KS2} have conjectured that, for 
$d\in\cal E$, the quantity $\log L(E_d,1)$ has a normal distribution with mean 
$-{1\over 2}\log \log|d|$ and variance $\log \log|d|$; see \cite{CKRS} \cite{djs} 
for numerical data towards this conjecture.  

Below we consider the families of quadratic twists $E_d$, where $E$ 
is one of the curve $A$, $B$, $C$, $D$, and $d$ runs over appropriate 
sets of positive integers dependent of $E$.  
Our data suggest that the 
values $\log L(E_d,1)$ also follow an approximate normal distribution. 
Let $W_E = \{d \leq 32\cdot 10^9 : d \, \, 
\text{satisfies ($**_E)$} \}$ 
and $I_x = [x,x+0.1)$ for $x\in\{ -10, -9.9, -9.8, \ldots, 10 \}$.
We create histograms with bins $I_x$ from the data 
$\left\{ \left(\log L(E_d,1) + \frac 1 2 \log\log d\right)/\sqrt{\log\log d} : d\in W_E\right\}$. 
Below we picture these histograms.

\begin{figure}[H]  
\centering
\includegraphics[trim = 0mm 20mm 0mm 15mm, clip, scale=0.4]{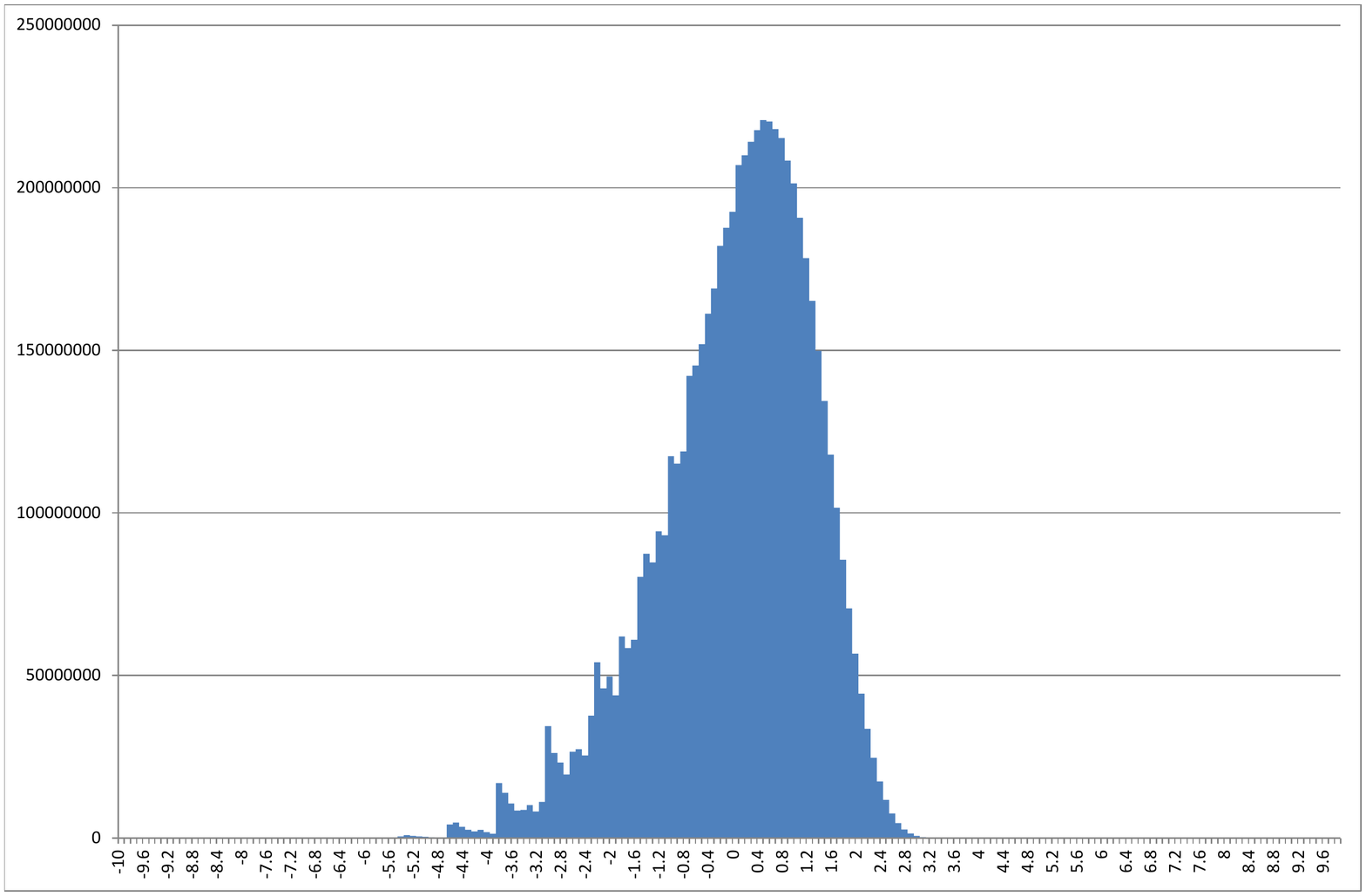} 
\caption{Histogram of values $\left( \log L(A_d,1) + {1\over 2} 
\log\log d \right)/\sqrt{\log\log d}$ for $d\in W_A$.} 
\end{figure}

\begin{figure}[H]  
\centering
\includegraphics[trim = 0mm 20mm 0mm 15mm, clip, scale=0.4]{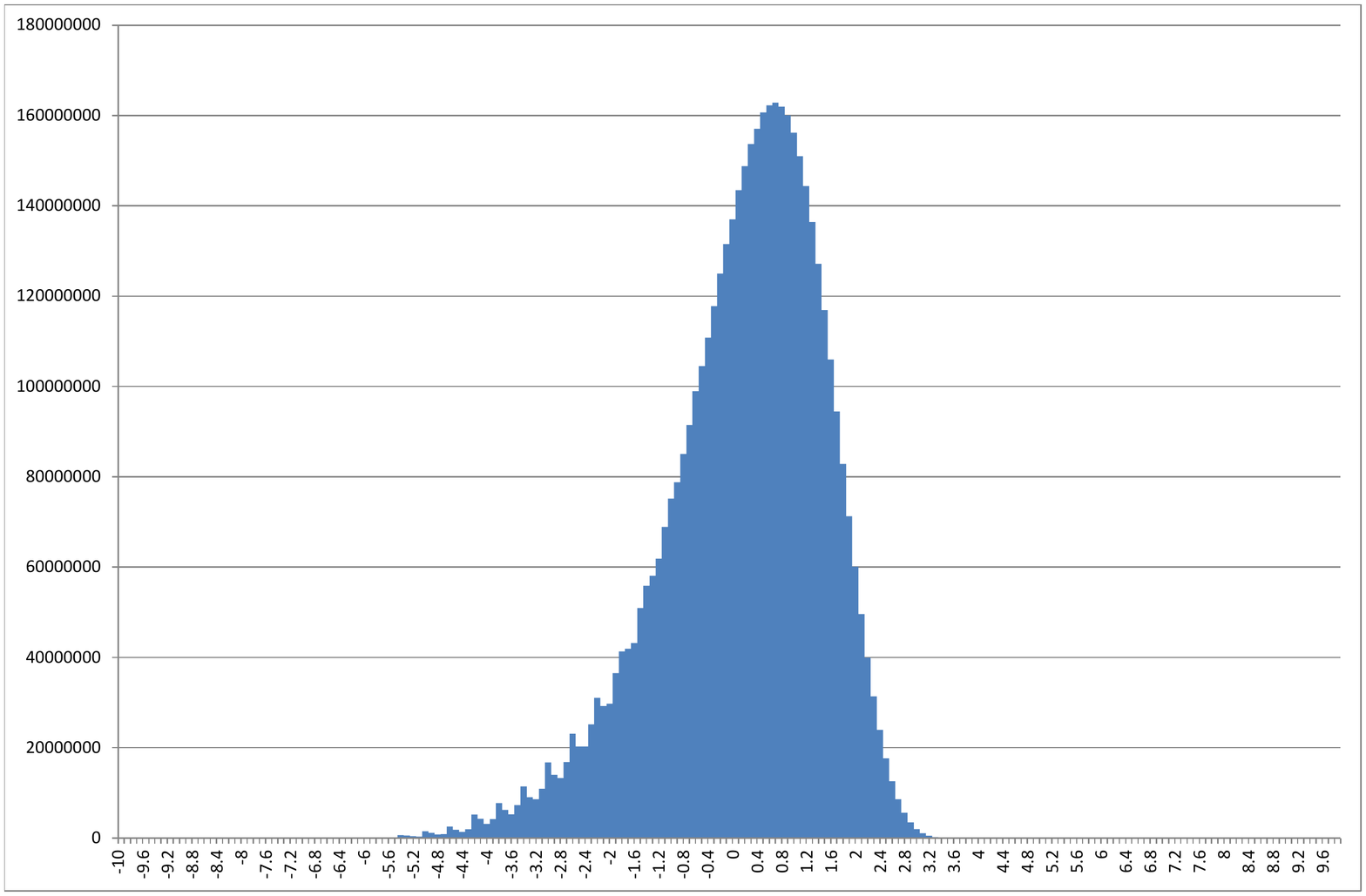} 
\caption{Histogram of values $\left( \log L(B_d,1) + {1\over 2} 
\log\log d \right)/\sqrt{\log\log d}$ for $d\in W_B$.} 
\end{figure}

\begin{figure}[H]  
\centering
\includegraphics[trim = 0mm 20mm 0mm 15mm, clip, scale=0.4]{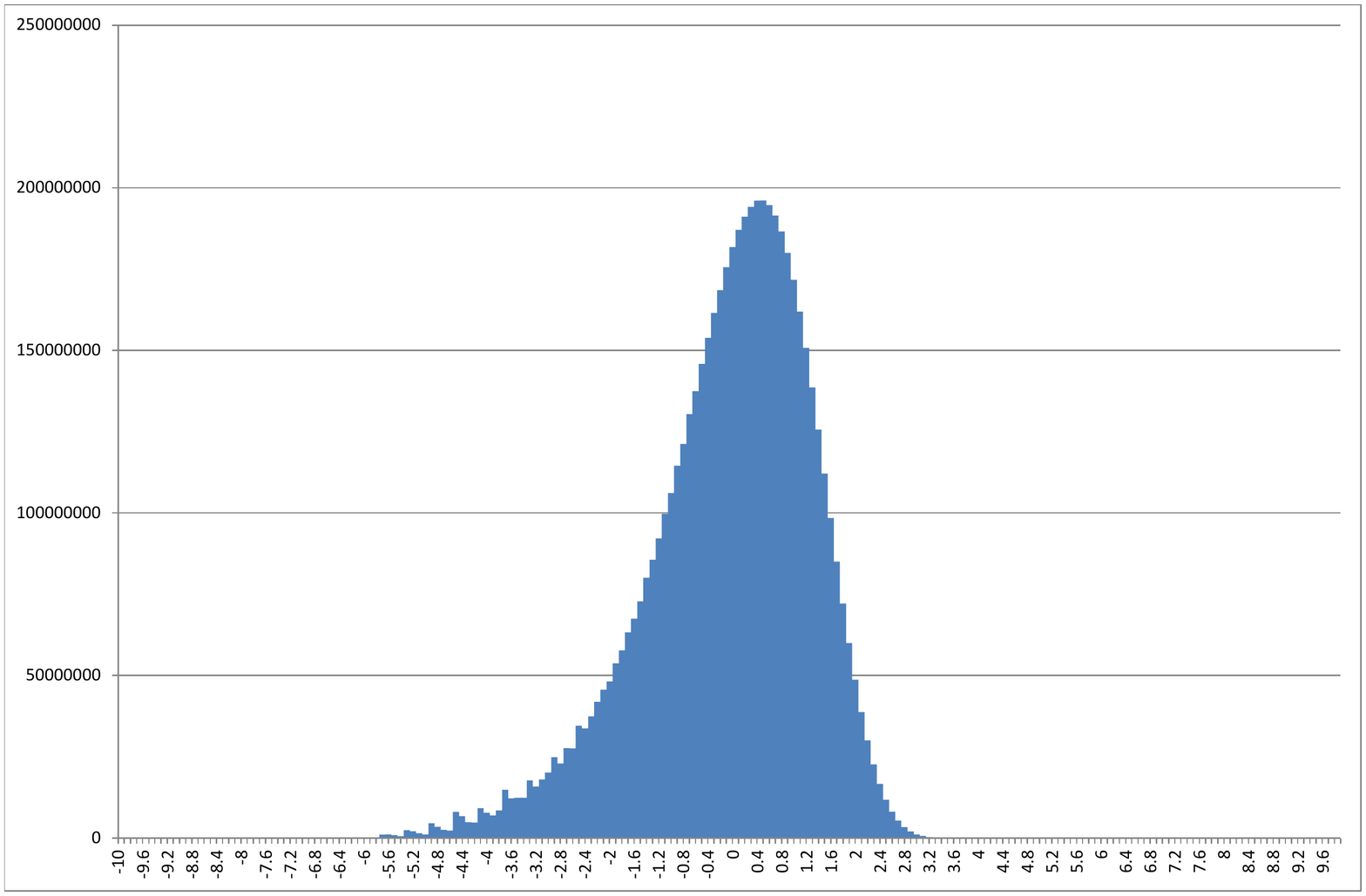} 
\caption{Histogram of values $\left( \log L(C_d,1) + {1\over 2} 
\log\log d \right)/\sqrt{\log\log d}$ for $d\in W_C$.} 
\end{figure}

\begin{figure}[H]  
\centering
\includegraphics[trim = 0mm 20mm 0mm 15mm, clip, scale=0.4]{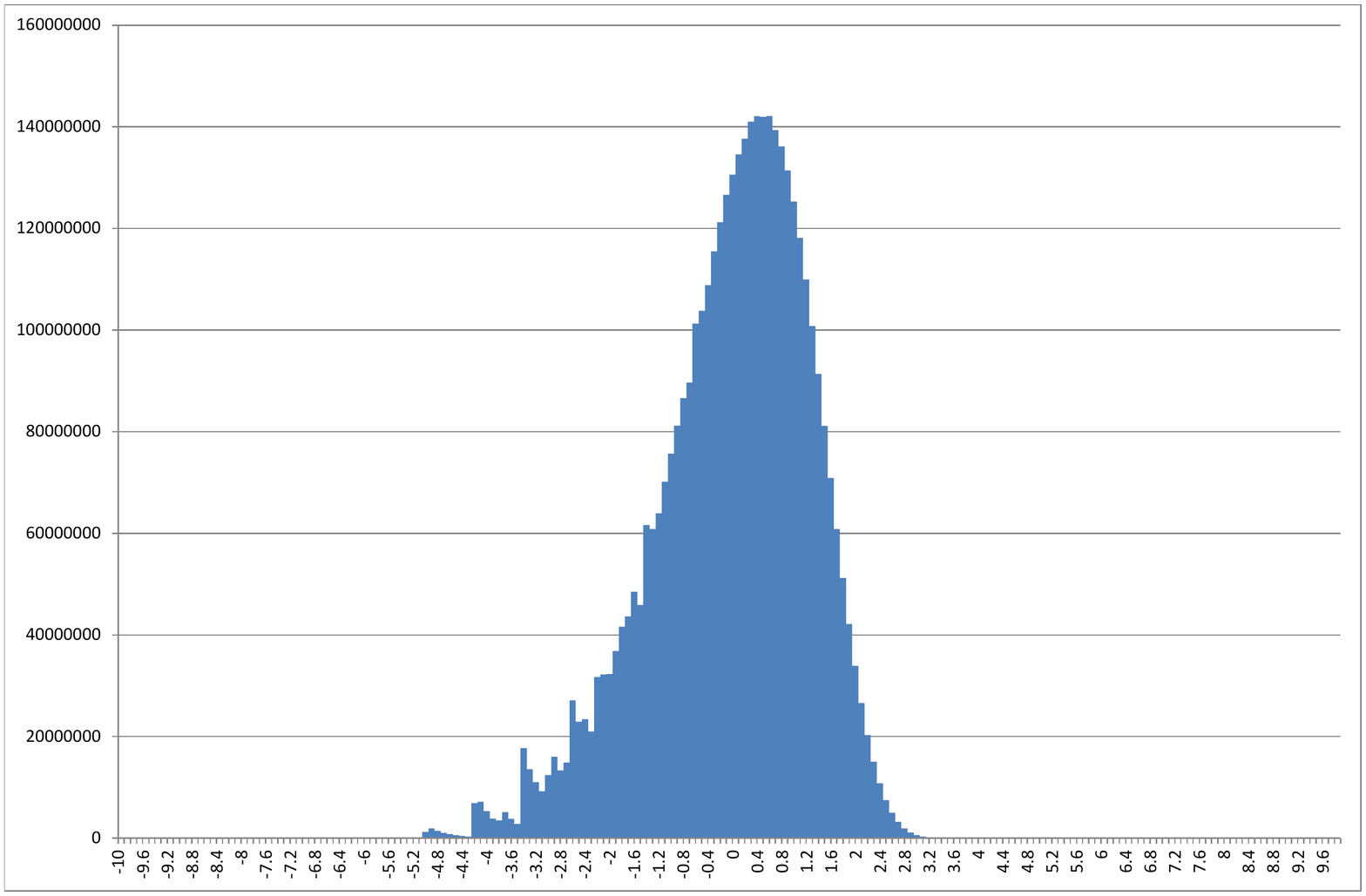} 
\caption{Histogram of values $\left( \log L(D_d,1) + {1\over 2} 
\log\log d \right)/\sqrt{\log\log d}$ for $d\in W_D$.} 
\end{figure}

\subsection{Distribution of $|\sza(E_d)|$} 

It is an interesting question to find 
results (or at least a conjecture) on distribution of the order of the Tate-Shafarevich 
group for rank zero quadratic twists of an elliptic curve over $\Bbb Q$. 
It turns out that the values of $\log(|\sza(E_d)|/\sqrt{d})$ 
are the natural ones to consider (compare Conjecture 1 in  \cite{RS}, 
and numerical experiments in \cite{djs}). 
Below we create histograms from the data 
$\left\{ \left( \log(|\sza(E_d)|/\sqrt{d}) - \mu_E  
\log\log d \right)/\sqrt{\sigma_E^2\log\log d} : \, d\in W_E\right\}$, 
for the curves $E\in\{A, B, C, D\}$, and 
where we take, according to Conjecture 1 in \cite{RS},  
$\mu_A=-{1\over 2}-2\log 2$, $\mu_B=\mu_D=-{1\over 2}-{3\over 2}\log 2$, 
$\mu_C=-{1\over 2}-{5\over 6}\log 2$,  
$\sigma_A^2=1+4(\log 2)^2$,  $\sigma_B^2=\sigma_D^2= 1+{5\over 2}(\log 2)^2$, 
and $\sigma_C^2=1+{7\over 6}(\log 2)^2$.

Our data suggest that the 
values $\log(|\sza(E_d)|/\sqrt{d})$ also follow an approximate normal 
distribution. Below we picture these histograms.

\begin{figure}[H]  
\centering
\includegraphics[trim = 0mm 20mm 0mm 15mm, clip, scale=0.4]{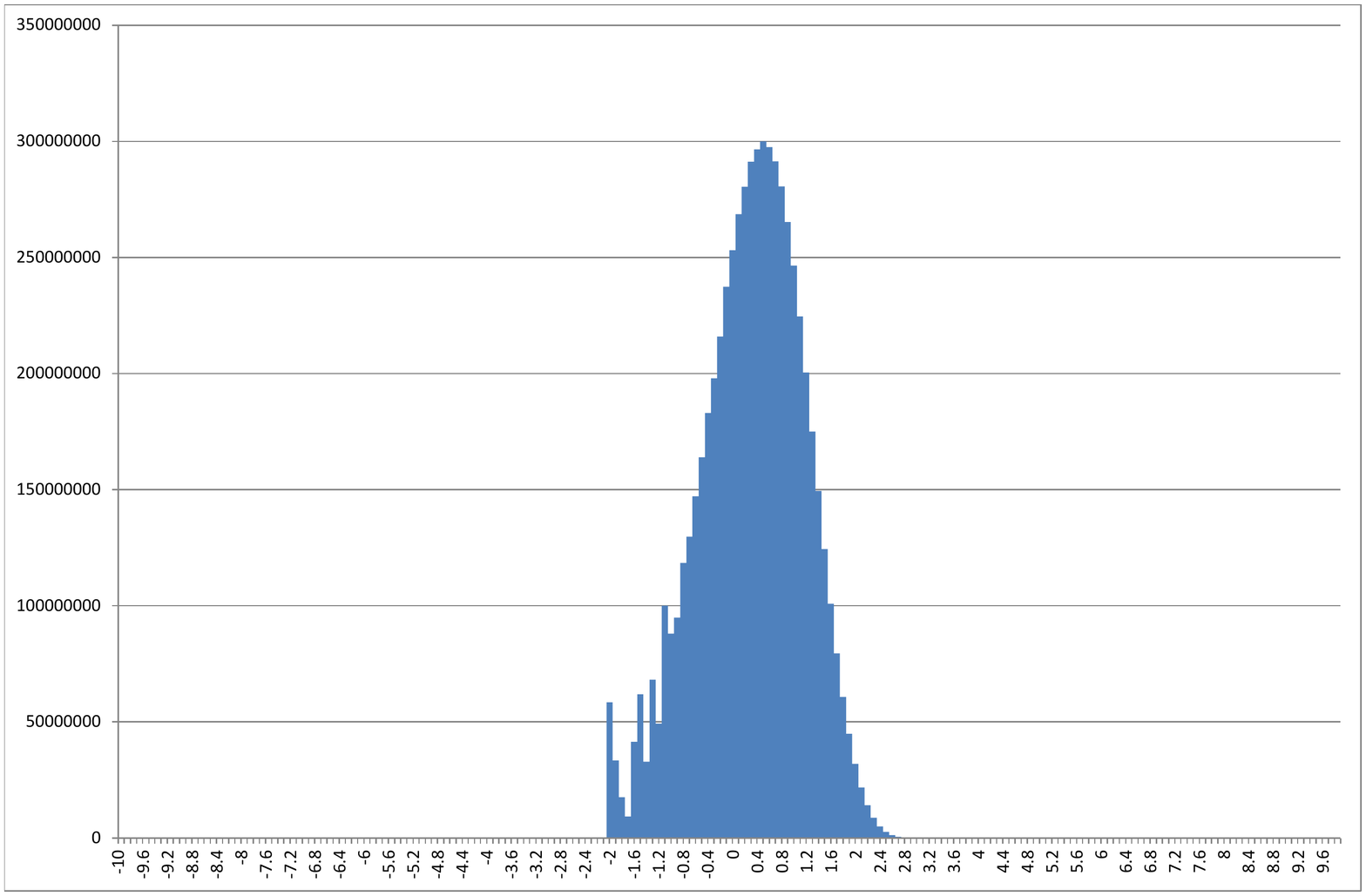} 
\caption{Histogram of values $\left( \log(|\sza(A_d)|/\sqrt d)  - \mu_A \log\log d \right)/\sqrt{\sigma_A^2 \log\log d}$ 
for $d\in W_A$.} 
\end{figure}

\begin{figure}[H]  
\centering
\includegraphics[trim = 0mm 20mm 0mm 15mm, clip, scale=0.4]{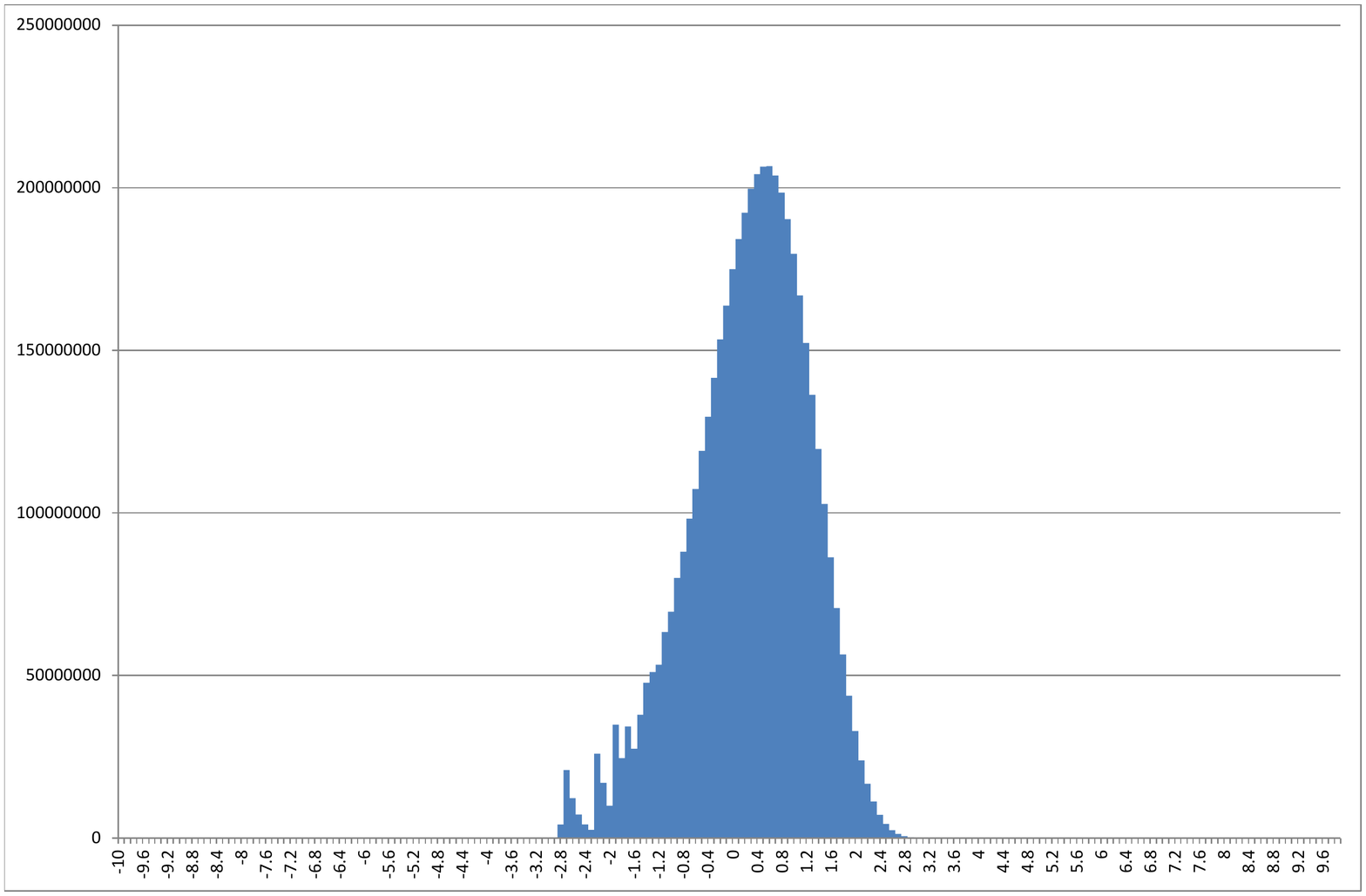} 
\caption{Histogram of values $\left( \log(|\sza(B_d)|/\sqrt d)  - \mu_B \log\log d \right)/\sqrt{\sigma_B^2 \log\log d}$ 
for $d\in W_B$.} 
\end{figure}

\begin{figure}[H]  
\centering
\includegraphics[trim = 0mm 20mm 0mm 15mm, clip, scale=0.4]{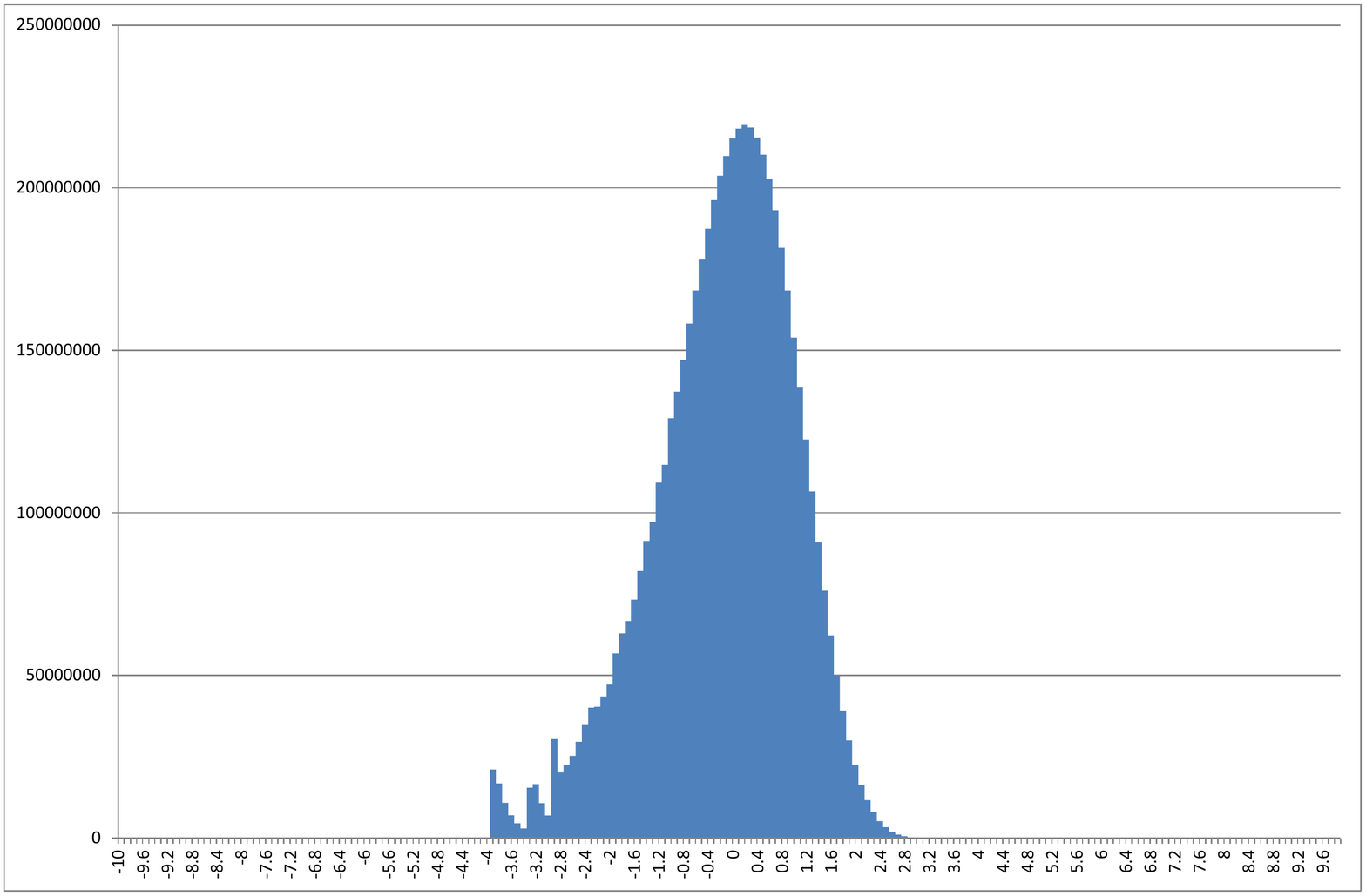} 
\caption{Histogram of values $\left( \log(|\sza(C_d)|/\sqrt d)  - \mu_C \log\log d \right)/\sqrt{\sigma_C^2 \log\log d}$ 
for $d\in W_C$.} 
\end{figure}

\begin{figure}[H]  
\centering
\includegraphics[trim = 0mm 20mm 0mm 15mm, clip, scale=0.4]{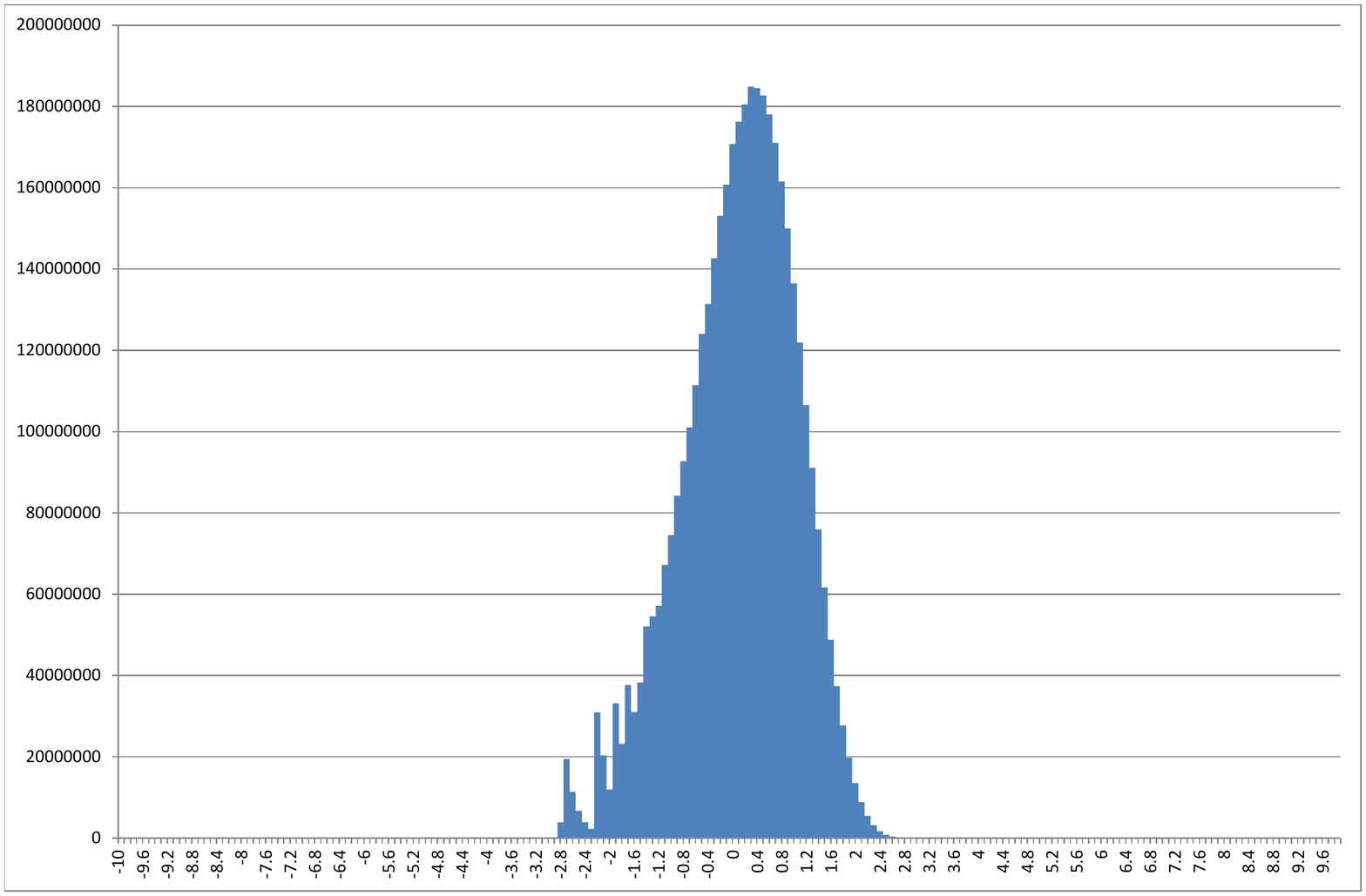} 
\caption{Histogram of values $\left( \log(|\sza(D_d)|/\sqrt d)  - \mu_D \log\log d \right)/\sqrt{\sigma_D^2 \log\log d}$ 
for $d\in W_D$.} 
\end{figure}

\section{More observations}

\subsection{Additional support towards the Conjecture 2}

Let $E\in\{A, B, C, D, X_0(49)\}$. 
For any elliptic curve $F$ in the isogeny 
class of $E$, we were able to confirm numerically the Conjecture 2. In all 
these cases we obtain $c_k(F) > 0$, for $1\leq k\leq 7$.

For any integer $l\geq 2$, let $d_l(E)$ denote 
the minimal $1\leq d \leq 32\cdot 10^9$, satisfying $(**_E)$, and such that 
$|\sza(E_d)|=l^2$. For each $2\leq l \leq 12$, we have considered quadratic twists 
of $E_{d_l(E)}$ by ${d\over d_l(E)}$, with $d\leq 32\cdot 10^9$ satisfying 
$(**_E)$. For any of these $55$ elliptic curves $E_{d_l(E)}$, we were able to 
confirm numerically the Conjecture 2 as well. In all these cases we obtain 
$c_k(E_{d_l(E)}) > 0$, for $1\leq k\leq 12$. 

One extra observation is that, in some cases, the number of twists with even 
order of $\sza$ is much larger than the number of twists with odd order of $\sza$. 
Examples include: $B_{d_8(B)}$, $D_{d_8(D)}$, and $E_{d_l(E)}$, where 
$E=X_0(49)$, $l=2, 4, 6, 12$.

\subsection{On a question of Coates-Li-Tian-Zhai}

Coates et al. \cite{cltz} have proved (among others) the following remarkable 
result. Let $E=X_0(49)$. 

\begin{thm} (\cite{cltz}, Theorem 1.2) 
Let $d=p_1\cdots p_l$ be a product of $\geq 0$ distinct primes, which 
are $\equiv 1$ mod $4$ and inert in $\Bbb Q(\sqrt{-7})$. Then 
$L(E_d,1)\neq 0$, $E_d(\Bbb Q)$ is finite, the Tate-Shafarevich 
group of $E_d$ is finite of odd cardinality, and the full 
Birch-Swinnerton-Dyer conjecture is valid for $E_d$.  
\end{thm} 

At the end of his paper, the authors say: ``{\it for every elliptic curve 
$E$ defined over $\Bbb Q$, we believe there should be some analogues of 
Theorems 1.1 - 1.4 for the family of quadratic twists of $E$, and it 
seems to us to be an important problem to first formulate precisely 
what such analogues should be, and then to prove them}.''

\bigskip 
Below we formulate analogues of Theorem 1.2 for the curves $A$ and $C$. 

\bigskip 

Let $E=A$. Let $d=p_1\cdots p_l$. Define the directed graph $G(d)$ by 
$V(G(d)):=\{p_1,...,p_l\}$ and 
$E(G(d)):=\{\overrightarrow{p_ip_j}: \, ({p_i\over p_j})=-1, \, 
1\leq i\not= j\leq l\}$. $G(d)$ is called odd provided that its only 
even partitions are trivial. Zhao \cite{Zhao} proved the 
following remarkable result. 

\begin{thm}
Suppose $d=p_1\cdots p_l$, $p_1\equiv 3 (\text{mod}\, 8)$, 
$p_2,...,p_l\equiv 1 (\text{mod}\, 8)$. Then $L(E_d,1)\neq 0$, 
$E_d(\Bbb Q)$ is finite, and the Tate-Shafarevich group of $E_d$ is 
finite of odd cardinality if $G(d)$ is odd. If $G(d)$ is odd, then 
the full Birch-Swinnerton-Dyer conjecture is valid for $E_d$. 
\end{thm}

\bigskip 
Let $E=C$, and $K=\Bbb Q(E[2])$.  Then Corollary 3.11 in \cite{Chen} 
allows to formulate the following 

\begin{thm}  
Let $d=p_1\cdots p_l$ be a product of $\geq 0$ distinct odd primes, which 
are prime to $11$ and split into two primes over $K$. Then 
$L(E_d,1)\neq 0$, $E_d(\Bbb Q)$ is finite, and the Tate-Shafarevich 
group of $E_d$ is finite of odd cardinality.  
\end{thm}

\subsection{Some folklore conjecture}

It is known \cite{Dab} that for any positive integer $m$, there are pairwise 
non-isogenous elliptic curves $E^1,\dots ,E^m$ defined over the rationals such 
that the rank of the Mordell-Weil group of the $p$-twist of $E^i$, $i=1,\dots,m$, 
has rank zero for a positive proportion of primes $p$. We propose the following 
(optimistic) 

\begin{conj}
For any positive integers $m$ and $k$, there are pairwise 
non-isogenous elliptic curves $E^1,\dots ,E^m$ defined over the rationals such 
that the rank of the Mordell-Weil group of the $p$-twist of $E^i$   
has rank zero for a positive proportion of primes $p$, and $|\sza(E^i)|=k^2$, 
for all $i=1,\dots,m$. 
\end{conj}

Our data for $E\in\{A,B,C,D\}$ (and for $E=X_0(49)$ in \cite{djs}) support 
the above conjecture for $m=3$ and small values of $k$.

\section{Elliptic curves with exceptionally large analytic order 
of the Tate-Shafarevich groups}

It has been known for a long time that $|\sza(E)|$ 
(provided is finite) can take arbitrarily large values (Cassels). 
The previously largest value for $|\sza(E)|$ was $63408^2$, found by 
D\k{a}browski and Wodzicki \cite{dw}. In (\cite{djs}, section 5) we 
propose a candidate with $|\sza(E)| > 100000^2$.  
Below we present the results of our search for 
elliptic curves with exceptionally large analytic order of the 
Tate-Shafarevich groups. We exibit $88$  examples of rank zero 
elliptic curves with $|\sza(E)| > 63408^2$. Our record is 
an elliptic curve $E=E_2(23,-348)$ with $|\sza(E)| = 1029212^2$. 
It is proven that it is the true order of $\sza$ in this case (see 8.4). 
Also  note that the prime $19861$ divides the orders of $\sza(E_i(22,304))$ - 
the largest (at the moment) prime dividing the order of $\sza(E)$ of 
an elliptic curve over $\Bbb Q$.

\subsection{Preliminaries} 

In this section we compute the analytic order of $\sza(E)$, i.e., the quantity 
$$
|\sza(E)|={L(E,1)\cdot |E(\Bbb Q)_{\text{tors}}|^2\over
C_{\infty}(E)C_{\text{fin}}(E)},
$$
for certain special curves of rank zero. We use the following approximation of
$L(E,1)$ 
$$
S_m=2\sum_{n=1}^m{a_n\over n}e^{-\frac{2\pi n}{\sqrt{N}}},
$$
which, for
$$
 m\geq {\sqrt{N}\over 2\pi}\left(2\log 2+k\log 10-\log(1-e^{-2\pi/\sqrt{N}})\right),
$$
differs from $L(E,1)$ by less than $10^{-k}$.

\smallskip
Consider (as in \cite{dw}) the family
$$
E_1(n,p):\quad y^2=x(x+p)(x+p-4\cdot 3^{2n+1}),
$$
with 
$(n,p)\in\Bbb N\times(\Bbb Z\setminus\{0\})$. Any member of the family admits
three isogenous (over $\Bbb Q$) curves $E_i(n,p)$ ($i=2,3,4$):

\bigskip 

$
E_2(n,p): \quad y^2=x^3+4(2\cdot 3^{2n+1}-p)x^2+16\cdot 3^{4n+2}x, 
$

$
E_3(n,p): \quad y^2=x^3+2(4\cdot 3^{2n+1}+p)x^2+(4\cdot 3^{2n+1}-p)^2x, 
$

$ 
E_4(n,p): \quad y^2=x^3+2(p-8\cdot 3^{2n+1})x^2+p^2x. 
$ 

\bigskip 
In our calculations, 
we focused on the pairs of integers $(n,p)$ within the bounds $20\leq n\leq 24$ 
and $0 < |p| \leq 5000$. Recall that the calculations in \cite{dw} were focused 
on the pairs $(n,p)$ within the bounds $3\leq n\leq 19$ and $0 < |p| \leq 1000$.

\bigskip 

The conductors, $L$-series and ranks of isogenous curves coincide, 
what may differ is the orders of $E(\Bbb Q)_{\text{tors}}$ and $\sza(E)$,
the real period $\Omega_E$, and the Tamagawa number $C_{\text{fin}}(E)$. 
In our situation we are dealing with $2$-isogenies, thus the  analytic order 
of $\sza(E)$ can only change by a power of $2$.

\bigskip 
{\it Notation}. Let $N(n,p)$ denote the conductor of the curve $E_i(n,p)$. 
We put $|\sza_i|=|\sza(E_i)|$.

\subsection{Elliptic curves $E_i(n,p)$ with 
$50000^2 \leq \max(|\sza_i|) < 250000^2$}

\begin{center}
\footnotesize
\begin{longtable}{|r|r|r|r|r|r|r|r|} 
\hline 
\multicolumn{1}{|c|}{$(n,p)$} 
& 
\multicolumn{1}{|c|}{$N(n,p)$} 
& 
\multicolumn{1}{|c|}{$|\sza_1|$} 
& 
\multicolumn{1}{|c|}{$|\sza_2|$} 
& 
\multicolumn{1}{|c|}{$|\sza_3|$} 
& 
\multicolumn{1}{|c|}{$|\sza_4|$} 
\\ 
\hline 
\endhead 
\hline 
\multicolumn{6}{|r|}{{Continued on next page}} \\
\hline
\endfoot

\hline 
\hline
\endlastfoot 

$(20,-756)$ & 42551829106699251024 & $27993^2$ &  $55986^2$ & $27993^2$  & $27993^2$  \\ 

$(20,-2000)$ & 190293894141760627320 & $15081^2$ & $60324^2$  &  $15081^2$ &  $60324^2$ \\

$(20,192)$ & 109418989131512359065 & $3780^2$ & $60480^2$  &  $945^2$ &  $60480^2$ \\ 

$(22,-692)$ & 11978814802342833513168 & $15194^2$ &  $30388^2$ &  $7597^2$ & $60776^2$  \\

$(21,-128)$  & 1969541804367222465954 &  $34234^2$  &  $68468^2$  & $34234^2$  & $68468^2$  \\ 

$(20,-180)$ & 60788327295284644080 & $20970^2$ & $41940^2$  & $10485^2$  &  $83880^2$ \\ 
 
$(21, 3)$  &   31512668869875559452120  &   $10962^2$  &  $43848^2$  &   $5481^2$  &   $87696^2$ \\ 

$(20,-2448)$ & 1653442502431742344680 &  $22028^2$ &  $88112^2$ &  $22028^2$ &  $88112^2$ \\ 

$(20,2704)$ & 11379574869677285146824 & $48538^2$ & $97076^2$  &  $97076^2$ &  $48538^2$ \\

$(21,12)$ & 281363114909603209392 & $12768^2$ & $102144^2$ & $3192^2$ &  $102144^2$ \\ 

$(20,-608)$ & 16631686347989878669080 & $25787^2$  & $103148^2$  & $51574^2$  & $51574^2$  \\

$(21,192)$ & 984770902183611232737 &  $54648^2$ & $109296^2$  &  $27324^2$ &  $109296^2$ \\ 

$(20, 4788)$  &  25871512096873143639456  &  $27745^2$  &  $110980^2$  &  $27745^2$  &  $27745^2$ \\ 

$(20,2680)$  &  23938195173261478962720 &  $14474^2$  &  $115792^2$  &  $14474^2$  &  $57896^2$  \\  

$(20,-801)$  &  34625031227394133415352  &  $29338^2$  &  $58676^2$  &  $29338^2$  &  $117352^2$  \\ 

$(22, 1344)$   &  62040566837567507664447  & $60930^2$  & $121860^2$ & $30465^2$ & $60930^2$   \\ 

$(20,-1436)$ & 1832369310703810488288 & $32455^2$ & $129820^2$  & $32455^2$  & $129820^2$  \\ 

$(20,4768)$ & 10032879618827902147272 & $16254^2$ & $130032^2$  & $8127^2$  & $65016^2$  \\  

$(21,-24)$  &  31512668869875559452768  &  $34092^2$  &    $68184^2$  &    $17046^2$  &  $136368^2$  \\  

$(20,-1376)$  &  37640132261240251922904  &  $70010^2$  &  $140020^2$  &  $140020^2$  &  $70010^2$  \\ 

$(22,64)$ & 8862938119652501095881 & $72306^2$ & $144612^2$  & $36153^2$  & $144612^2$  \\ 

$(21,-1536)$ & 1969541804367222468066 & $75897^2$ &  $151794^2$ & $75897^2$  &  $151794^2$ \\ 
 
$(20,-6)$ & 14005630608833581979328 & $19248^2$ & $76992^2$  & $4812^2$  &  $153984^2$ \\ 

$(22,304)$  &  27493195799738370745848  &  $39722^2$  &  $158888^2$  &    $19861^2$  &  $79444^2$  \\  

$(21,1516)$ & 11663380372737145525968 & $25866^2$ &  $206928^2$ &  $12933^2$ &  $103464^2$ \\ 

$(21,480)$  &  39390836087344449300840  &  $54110^2$  &  $216440^2$  &  $27055^2$  &  $108220^2$  \\ 

$(23,1452)$ & 6451697601805864768272 & $55698^2$ & $222792^2$  & $27849^2$  &  $222792^2$ \\

\end{longtable}
\end{center}

\subsection{Elliptic curves $E_i(n,p)$ with 
$\max(|\sza_i|) \geq 250000^2$}

\begin{center}
\footnotesize
\begin{longtable}{|r|r|r|r|r|r|r|r|} 
\hline 
\multicolumn{1}{|c|}{$(n,p)$} 
& 
\multicolumn{1}{|c|}{$N(n,p)$} 
& 
\multicolumn{1}{|c|}{$|\sza_1|$} 
& 
\multicolumn{1}{|c|}{$|\sza_2|$} 
& 
\multicolumn{1}{|c|}{$|\sza_3|$} 
& 
\multicolumn{1}{|c|}{$|\sza_4|$} 
\\ 
\hline 
\endhead 
\hline 
\multicolumn{6}{|r|}{{Continued on next page}} \\
\hline
\endfoot

\hline 
\hline
\endlastfoot

$(21,4)$ & 15756334434937779726048 & $130614^2$ &  $261228^2$ & $65307^2$  & $261228^2$  \\ 

$(21,1248)$    &  102416173827095568122280  &  $70375^2$  &  $281500^2$  &  $70375^2$  &  $70375^2$  \\  

$(20,-201)$    &  234594312697962498467304 &  $141540^2$  &  $141540^2$  &  $283080^2$  &  $141540^2$   \\

$(23,960)$ &  398832215384362549313205 &  $96254^2$ & $385016^2$  & $48127^2$  & $192508^2$  \\ 

$(24,832)$ & 373306953599763346160205 & $75780^2$ & $303120^2$ & $37890^2$ & $151560^2$  \\   

$(20,1120)$ & 30637316956823460343320 & $20440^2$ &  $327040^2$ & $20440^2$  &  $81760^2$ \\ 

$(23,-84)$ & 17448909423065861532624 & $184991^2$ & $369982^2$  & $184991^2$  &  $184991^2$ \\

$(22,480)$  &  354517524786100043822760  & $99938^2$   & $399752^2$   & $49969^2$   & $199876^2$ \\ 

$(23,-8)$ & 7441767284139709375008  &  $102120^2$ & $204240^2$  & $51060^2$  &  $408480^2$ \\ 

$(21,-233)$   &  149845956054714394972728  &  $51317^2$  &  $205268^2$  &  $51317^2$ & $410536^2$ \\ 

$(23,-96)$  &  638131544614980078907464   &   $264696^2$  &  $529392^2$  &  $132348^2$  &  $529392^2$  \\ 

$(24,-96)$ & 302272836922885300534872  & $412146^2$  &  $824292^2$ &  $206073^2$ & $824292^2$  \\  

$(23,-348)$ & 37011629587668844576720608 & $514606^2$ & $1029212^2$ & $257303^2$ & $1029212^2$  \\ 

\end{longtable}
\end{center}

\subsection{Birch and Swinnerton-Dyer conjecture for elliptic curves with 
exceptionally large analytic order of Tate-Shafarevich groups}

In this subsection, we will use the deep results by Skinner-Urban 
\cite{SkUr}, to prove the full 
version of the Birch-Swinnerton-Dyer conjecture for some  
elliptic curves $E_i(n,p)$ with exceptionally large analytic order 
of Tate-Shafarevich groups.  

Let $\overline{\rho}_{E,p}: 
\text{Gal}(\overline{\Bbb Q}/\Bbb Q) \to 
\text{GL}_2(\Bbb F_p)$ denote the Galois 
representation on the $p$-torsion of $E$. 
Assume $p\geq 3$. 

\begin{thm} (\cite{SkUr}, Theorem 2) Let $E$ be an 
elliptic curve over $\Bbb Q$ with conductor $N_E$. 
Suppose: (i) $E$ has good ordinary reduction at $p$; 
(ii) $\overline{\rho}_{E,p}$ is irreducible; (iii) 
there exists a prime $q\not=p$ such that $q\mid\mid N_E$ 
and $\overline{\rho}_{E,p}$ is ramified at $q$; 
(iv) $\overline{\rho}_{E,p}$ is surjective.  
If moreover $L(E,1)\not=0$, then the $p$-part of the Birch 
and Swinnerton-Dyer conjecture holds true, and  
we have 
$$
\text{ord}_p(|\sza(E)|)=\text{ord}_p\left(
{|E(\Bbb Q)_{\text{tors}}|^2 L(E,1)\over C_{\infty}(E)C_{fin}(E)}\right). 
$$
\end{thm}

The following curves from the tables in 8.2 and 8.3 satisfy the conditions 
$(N_{E_i},\min (|\sza_i|))=1$:  
$E_i(22,692)$, $E_i(20,-608)$, $E_i(20,-1436)$, $E_i(23,-84)$, $E_i(20,4788)$, 
$E_i(22,304)$, $E_i(21,1248)$, $E_i(22,480)$, $E_i(23,960)$ and $E_i(23,-348)$. 
In these cases, all the assumptions of the result by Skinner-Urban are 
satisfied. 

Let us give some details for the curves  $E_i=E_i(20,-1436)$. 
We can use Theorem 2 
to show that $|\sza_1|=5^26491^2$ is the true order of 
$\sza(E_1)$ (and, hence, all $|\sza_i|$ are the true orders 
of $\sza(E_i)$). 
(i) $E_1$ has good ordinary reduction at $5$ 
and $6491$: $(N_{E_1},5)=(N_{E_1},6491)=1$, and $a_5(E_1)=2$, 
$a_{6491}(E_1)=108$. (ii) We have the following general result of 
B. Mazur. Let $E$ be an elliptic curve over 
$\Bbb Q$ with all its $2$-division points defined over $\Bbb Q$. 
Then $\overline{\rho}_{E,p}$ is absolutely irreducible for any 
prime $p\geq 5$. (iii) Take $q=7$. Then $7 || N_{E_1}$, and 
$\overline{\rho}_{E_1,p}$ is ramified at $7$ for $p=5$ and $6491$, 
since these $p$'s do not divide $\text{ord}_7(\Delta_{E_1})$. 
(iv) We have  
$j_{E_1}={2^6 285451^3 4660272567723053424015171049538317^3 
\over 3^{82} 7^8 31^2 257^2 359^2 107323^2 8883041^2}$. 
The representation  $\overline{\rho}_{E_1,p}$ is surjective 
for any prime $p\geq 19$ by Prop. 1.8 in \cite{Zyw}. 
On the other hand, Prop. 6.1 in \cite{Zyw}  
gives a criterion to determine whether $\overline{\rho}_{E,p}$ 
is surjective or not for any non-CM elliptic curve $E$ 
and any prime $p\leq 11$. For instance, the representation 
$\overline{\rho}_{E,5}$ is not surjective if and only if 
$j_E={5^3(t+1)(2t+1)^3(2t^2-3t+3)^3 \over (t^2+t-1)^5}$ or 
$j_E={5^2(t^2+10t+5)^3\over t^5}$ or $j_E=t^3(t^2+5t+40)$ 
for some $t\in\Bbb Q$. Taking 
$t={a\over b}$, where $a$, $b$ are relatively prime integers 
with  $b>0$, we obtain 
$$
{5^3(a+b)(2a+b)^3(2a^2-3ab+3b^2)^3 \over (a^2+ab-b^2)^5} = 
j_{E_1}  \quad \text{or}  
$$
$$
{5^2(a^2+10ab+5b^2)^3 \over a^5b} = j_{E_1} \quad \text{or} 
\quad {a^3(a^2+5ab+40b^2) \over b^5} = j_{E_1}. 
$$
In the first case note, that $a^2+ab-b^2$ is relatively prime 
to $3$, a contradiction. The second case is impossible, since 
necessarily $ab\equiv 1 (\text{mod}\,8$, and consequently $8$ divides 
$a^2+10ab+5b^2$, a contradiction with $\text{ord}_2(j_{E_1})=6$.  
The last case is impossible, since a denominator of $j_{E_1}$ 
is not a fifth power of an integer.

\section{Examples of elliptic curves with $c_{2m+1}(E)=0$ and $c_{4m+2}(E)=0$}

\subsection{Numerical data}

Let $E^i=E_i(5,2)$. For all positive, square-free integers $d\leq 2495$, 
prime to $N_{E^i}$, and such that the groups $E^i_d(\Bbb Q)$ are finite, 
we calculated the orders of $\sza(E^i_d)$. A part of the data is collected 
below in the following table.

\begin{center}
\footnotesize
\begin{longtable}{|r|r|r|r|r|r|r|} 
\hline 
\multicolumn{1}{|c|}{$d$} 
& 
\multicolumn{1}{|c|}{$|\sza(E^1_d)|$} 
& 
\multicolumn{1}{|c|}{$|\sza(E^2_d)|$} 
& 
\multicolumn{1}{|c|}{$|\sza(E^3_d)|$} 
& 
\multicolumn{1}{|c|}{$|\sza(E^4_d)|$} 
\\ 
\hline 
\endhead 
\hline 
\multicolumn{5}{|r|}{{Continued on next page}} \\
\hline
\endfoot

\hline 
\hline
\endlastfoot 
 
$1$  & $3^2$  & $6^2$  & $3^2$  & $12^2$  \\ 
$5$  &  $6^2$ & $24^2$ &  $3^2$  & $24^2$ \\ 
$11$  &  $6^2$  & $24^2$  &  $3^2$  & $24^2$ \\ 
$23$  &  $6^2$  & $24^2$  &  $6^2$  & $24^2$ \\ 
$55$  &  $3^2$  & $12^2$  &  $3^2$  & $24^2$ \\ 
$59$  &  $12^2$  & $48^2$  &  $6^2$  & $48^2$ \\ 
$61$  &  $6^2$  & $12^2$  &  $3^2$  & $24^2$  \\ 
$71$  &  $6^2$  & $24^2$  &  $3^2$  & $24^2$  \\ 
$73$  &  $2^2$  &  $4^2$  &  $1^2$  &  $8^2$  \\ 
$79$  &  $4^2$  &  $8^2$  &  $2^2$  & $16^2$  \\ 
$83$  &  $2^2$  &  $8^2$  &  $2^2$  &  $8^2$  \\ 
$97$  &  $2^2$  &  $4^2$  &  $1^2$  &  $8^2$  \\ 
$101$  &  $20^2$  & $80^2$  & $10^2$  & $80^2$  \\ 
$109$  &  $4^2$  &  $8^2$  &  $2^2$  & $16^2$  \\ 
$113$  &  $8^2$  & $32^2$  &  $4^2$  & $32^2$  \\ 
$115$  &  $2^2$  &  $8^2$  &  $1^2$  &  $8^2$  \\ 
$119$  &  $8^2$  & $32^2$  &  $8^2$  & $64^2$  \\ 
$127$  &  $8^2$  & $16^2$  &  $4^2$  & $32^2$  \\ 
$143$  &  $3^2$  & $12^2$  &  $3^2$  & $24^2$  \\ 
$157$  &  $2^2$  &  $4^2$  &  $1^2$  & $8^2$  \\ 
$163$  &  $6^2$  &  $12^2$  &  $3^2$  & $24^2$  \\ 
$173$  &  $2^2$  &  $8^2$  &  $1^2$  &  $8^2$ \\ 
$179$  &  $28^2$  & $112^2$  & $14^2$ & $112^2$  \\ 
$197$  &  $10^2$  & $40^2$  &  $5^2$  & $40^2$  \\ 
$199$  &  $1^2$  &  $2^2$  &  $1^2$  &  $4^2$  \\ 
$205$  &  $4^2$  & $16^2$  &  $2^2$  & $16^2$  \\ 
$217$  &  $2^2$  & $8^2$  &  $2^2$ &  $16^2$  \\ 
$227$  &  $4^2$  & $16^2$  &  $2^2$  & $16^2$  \\ 
$235$  &  $2^2$  &  $8^2$  &  $1^2$  &  $8^2$  \\ 
$253$  &  $2^2$  &  $8^2$  &  $1^2$  &  $8^2$  \\ 
$257$  &  $1^2$  &  $4^2$  &  $1^2$  &  $4^2$  \\ 
$263$  &  $6^2$  & $24^2$  &  $3^2$  & $24^2$  \\ 
$271$  &  $3^2$  &  $6^2$  &  $3^2$  & $12^2$  \\ 
$277$  &  $9^2$  & $18^2$  &  $9^2$  & $36^2$  \\ 
$281$  &  $27^2$ & $108^2$  & $27^2$ & $108^2$  \\ 
$283$  &  $2^2$  &  $4^2$  &  $1^2$  &  $8^2$  \\ 
$295$  &  $3^2$  & $12^2$  &  $3^2$  & $24^2$  \\ 
$299$  &  $2^2$  &  $8^2$  &  $1^2$  &  $8^2$  \\ 
$301$  &  $2^2$  &  $8^2$  &  $1^2$  &  $8^2$  \\ 
$305$  &  $1^2$  &  $4^2$  &  $1^2$  &  $8^2$  \\  

\end{longtable}
\end{center} 

Our data strongly suggest that 
$c_{2m+1}(E^2_d) = c_{2m+1}(E^4_d) = c_{4m+2}(E^2_d) = 0$. 
We prove these statements in the next subsection.

\subsection{Lower bounds for the $2$-rank of $\sza(E'_d)$}

Consider an elliptic curve $E=E_1(5,2)$ given by the equation 
$y^2=x^3 - 708584x^2 - 1417172x$, 
of conductor $N_E=2^6.3.19.29.643$.  
Let $E'=E_2(5,2)$ be the elliptic curve given by the equation  
$y^2=x^3 + 1417168x^2+502096953744x$. 
 Let $E_d$ and $E'_d$ denote the quadratic twists by positive, square-free 
integers $d$, prime to $N_E$. Consider the two-isogeny  
$\phi: E_d\to E'_d$, defined by 
$\phi((x,y))=(y^2/x^2,-y(1417172r^2+x^2)/x^2)$; 
let $\hat\phi$ denote the dual isogeny. Consider the Selmer groups 
$S^{(\phi)}(E_d/\mathbf Q)$, and  
$S^{(\hat\phi)}(E'_d/\mathbf Q)$.  
We use the notations and results from  chapter X 
of Silverman's book \cite{sil}. Let 

$
C^{(d)}_r: ry^2=r^2 + 2^4.23.3851rdx^2 + 2^4.3^{22}d^2x^4, 
$

$
C^{'(d)}_r: ry^2=r^2 - 2^5.23.3851rdx^2 -2^6.19.29.643d^2x^4 
$

\noindent 
be the principal homogeneous spaces under the actions of the elliptic 
curves previously defined. Let $\Sigma(M)$ and $\Delta(M)$ be the  support 
of an integer $M$ in the set of prime numbers and the set 
of divisors of $M$ in $\mathbf Z$ respectively. Using 
[\cite{sil}, Proposition 4.9, p.302],  we have the following  identifications: 

\noindent 
$
S^{(\phi)}(E_d/\mathbf Q)\simeq \{r\in \Delta(2.3.d): 
C^{(d)}_r(\mathbf Q_l)\not=\emptyset \quad \forall \, 
l\in \Sigma(2.3.23.3851.d)\cup\{\infty\}\}, 
$

\noindent 
$
S^{(\hat\phi)}(E'_d/\mathbf Q) \simeq \{r\in \Delta(2.19.29.643.d): 
C^{'(d)}_r(\mathbf Q_l)\not=\emptyset \quad \forall \, 
l\in \Sigma(2.3.19.29.643.d)\cup\{\infty\}\}. 
$ 

\begin{thm}
Let $d$ be a positive, square-free integer, prime to $N_E$, and such 
that the group $E(\Bbb Q)$ is finite. Then the group $\sza(E'_d)[2]$ 
is non-trivial. 
\end{thm} 

{\it Proof}. We have  
$\dim_2 S^{(\phi)}(E_d/\mathbf Q)\geq \dim_2 \sza(E_d)[\phi]$. 
If $E_d(\Bbb Q)$ is finite, then using the fundamental formula 
(see, for instance, \cite{sil}, p. 314), we obtain 
$\dim_2 \sza(E'_d)[\hat\phi] \geq \dim_2 S^{(\hat\phi)}(E'_d/\mathbf Q) - 2$. 
Here $\sza(E_d)[\phi]$ is the kernel of the mapping $\sza(E_d) \to \sza(E'_d)$ 
induced by $\phi$, and $\sza(E'_d)[\hat\phi]$ is defined similarly; also, 
we write $\dim_2$ for $\dim_{\Bbb F_2}$. 

Now, it is sufficient to prove the following result. 

\begin{lem} 
We have $\dim_2 S^{(\hat\phi)}(E'_d/\mathbf Q) \geq 3$. 
\end{lem}

{\it Proof of Lemma 1}. We will exhibit $8$ elements in 
$S^{(\hat\phi)}(E'_d/\mathbf Q)$. More precisely, it is sufficient to check that 
$<19.643, -19.29,  r_d> \subset S^{(\hat\phi)}(E'_d/\mathbf Q)$, 
where $r$ denotes any prime divisor of $d$, and  

$$ 
r_d=
\begin{cases}
r,  \quad  \text{if} \quad  r\equiv 1(\text{mod}\, 8),\\ 
19r  \quad \text{if} \quad  r\equiv 3(\text{mod}\, 8),\\
29r  \quad \text{if} \quad  r\equiv 5(\text{mod}\, 8),\\ 
-r   \quad \text{if} \quad  r\equiv 7(\text{mod}\, 8).   
\end{cases}
$$ 

We omit the standard calculations using the Hensel Lemma. 

Similarly, one can prove that $c_{2m+1}(E^4_d)=0$.

Institute of Mathematics, University of Szczecin, Wielkopolska 15, 
70-451 Szczecin, Poland; E-mail addresses: dabrowsk@wmf.univ.szczecin.pl and dabrowskiandrzej7@gmail.com;  
lucjansz@gmail.com

\end{document}